\documentclass[12pt]{amsart}

\usepackage{amssymb,stmaryrd}
\usepackage{booktabs}
\usepackage{cases}
\usepackage{array} 
\usepackage{float} 
\usepackage{rotating }

\usepackage{indentfirst}
\usepackage{multirow}
\usepackage{amsthm}
\usepackage{longtable}
\usepackage{float}
\usepackage{arydshln}
\usepackage{mathtools}
\usepackage{multirow}

\usepackage[foot]{amsaddr}
\makeatletter
\renewcommand{\email}[2][]{%
  \ifx\emails\@empty\relax\else{\g@addto@macro\emails{,\space}}\fi%
  \@ifnotempty{#1}{\g@addto@macro\emails{\textrm{(#1)}\space}}%
  \g@addto@macro\emails{#2}%
}
\makeatother
\usepackage{mathtools}
\usepackage{enumerate,pdflscape}
\usepackage{multirow}
\usepackage{color}
  \usepackage[titletoc]{appendix}

\makeatletter
\newcommand{\svast}{\bBigg@{3}}
\newcommand{\vast}{\bBigg@{3.5}}
\newcommand{\Vast}{\bBigg@{5}}
\makeatother

 \newcommand{\ct}{C_2}
 \newcommand{\at}{A_2}
 \newcommand{\aoao}{A_1\times A_1}

 \newcommand{\s}{\mathfrak{sl}(2,\mathbb{C})}

 \newcommand{\ad}{\text{ad}}

  \theoremstyle{definition}

  \theoremstyle{plain}

\newcommand\borel{\mathfrak{b}}

\newcommand\algn{\mathfrak{n}}

\newcommand\algt{ \frak{h}}

\newcommand\ga{G_2}

\title[Subalgebras of semisimple Lie algebras]{Subalgebras of the rank two  semisimple Lie algebras}

\begin{document}
\date{\today}                                           

\author[Andrew Douglas]{Andrew Douglas$^{1,2,3}$}
\address[]{$^1$Department of Mathematics, New York City College of Technology, City University of New York, Brooklyn, NY, USA.}
\address[]{$^2$Ph.D. Programs in Mathematics and Physics, CUNY Graduate Center, City University of New York, New York, NY, USA.}

\author[Joe Repka]{Joe Repka$^3$}
\address{$^3$Department of Mathematics, University of Toronto, Toronto, ON, M5S 2E4, Canada}
\email{adouglas@citytech.cuny.edu, repka@math.toronto.edu}




\begin{abstract}
In this expository article, we describe the classification of the subalgebras of the rank $2$ semisimple Lie algebras. Their semisimple  subalgebras are well-known, and in a recent series of papers, we completed the classification
of the subalgebras
of the \textit{classical} rank 2 semisimple Lie algebras. Finally, Mayanskiy finished the classification of the subalgebras of the remaining rank $2$ semisimple Lie algebra, the exceptional Lie algebra $\ga$.  We identify subalgebras of the classification in terms of a uniform classification scheme of Lie algebras of low dimension. The classification is up to inner automorphism, and the ground field is the complex numbers.
\end{abstract}

\maketitle
\tableofcontents
\listoftables

\section{Introduction}

In this expository article, we describe the classification of the subalgebras of the rank $2$ semisimple Lie algebras. Specifically, we describe the subalgebras of the symplectic algebra $C_2$, the special linear algebra $A_2$, the semisimple--and not simple--Lie algebra $A_1\times A_1$, and the exceptional Lie algebra
$G_2$. The classification is up to inner automorphism, and the ground field is the complex numbers.

We omit virtually all proofs in this article, as they may be found in the literature, which we cite. But, we do outline techniques, strategy, and theory employed in proofs for one case--the subalgebras of $C_2$.

By Levi's Theorem [\cite{levi}, Chapter II, Section 2],  a subalgebra of a complex semisimple Lie algebra is either semisimple, solvable, or a nontrivial semidirect sum of the first two.  A subalgebra that is a  nontrivial semidirect sum of a semisimple subalgebra with a solvable subalgebra is called a Levi decomposable subalgebra. 

Semisimple subalgebras of semisimple Lie algebras have been extensively studied  \cite{degraafd, dynkin, dynkin2, lorent, min}.   For instance, the semisimple subalgebras of the exceptional Lie algebras have been classified, up to inner automorphism \cite{min}.  As another important example,   de Graaf \cite{degraafd} classified the semisimple subalgebras of the simple Lie algebras of ranks $\le 8$. In particular, the semisimple subalgebras of the rank $2$ simple Lie algebras have been known for some time.

Until quite recently, much less research had examined and classified the {\it non-semisimple} subalgebras of semisimple Lie algebras.  
We recently classified the solvable and Levi decomposable subalgebras of the symplectic algebra $C_2$ \cite{dc2, c2} and the special Linear algebra $A_2$ \cite{a2}.  We also gave a full classification of the subalgebras of the semisimple Lie algebra $A_1\times A_1$ \cite{drsof}. In $2016$, 
Mayanskiy posted a classification of the subalgebras of $G_2$ \cite{may}. The aforementioned work completes the classification of the subalgebras of the rank $2$ semisimple Lie algebras and a  description of this classification, and identification thereof, encompasses the present article. 

Interestingly, in Mayanskiy's \cite{may} classification of $G_2$ subalgebras, there are  naturally parameterized families of inequivalent subalgebras  for which subalgebras of the family are not only inequivalent 
(under inner automorphisms),
but are also non-isomorphic. This 
phenomenon
also occurs in 
the
classification of $A_2$ and $C_2$ subalgebras
\cite{a2, c2}, where 
we precisely describe the different isomorphism classes within such a family.
In this article, we 
also
identify such families for $G_2$, and the non-isomorphic classes therein.

When possible, we identify solvable and Levi decomposable subalgebras in terms of the (partial) classification of solvable and Levi decomposable algebras presented by \v{S}nobl and Winternitz \cite{levi}. Further, we describe the pertinent  portion of  \v{S}nobl and Winternitz's classification  in the appendices of the present article. 

Our rationale for writing this article is threefold. First, we felt it beneficial to students and 
researchers
in the field to have  the entire subalgebra classification in one well-organized article. Second, we wanted to describe the classification using one consistent choice of notation, which is not done at present in the literature on subalgebras. Thirdly, we wanted to identify the subalgebras with respect to one classification, which was not previously done in totality. Specifically, identification was not 
available for the subalgebras of $G_2$. 

We end the section by briefly articulating the importance of the mathematical domain of this article--subalgebras of Lie algebras. In addition to the intrinsic mathematical significance of  classifications of  subalgebras of Lie algebras (or classifications of subgroups of corresponding  Lie groups), such classifications also have physical significance and mathematical applications, some of which are listed below (as we described in \cite{c2}):
\begin{itemize}
\item If a system of differential equations is invariant under a Lie group, then its
subgroups can be used to construct group invariant solutions \cite{pjoliver}. 
\item Subgroups of the symmetry groups of nonlinear partial differential equations provide a method for performing
symmetry reduction (reducing  the number of independent
variables) \cite{wlong, ghp,dklw}.
\item A knowledge of the subgroup structure of a Lie group $G$ is needed if we are interested in considering all possible contractions of $G$ to other groups \cite{pw77}.
\item Physical models--such as the vibron model, and the interacting boson model-- use chains of subalgebras, and these subalgebras need to be explicitly described in 
applications
\cite{iachello}.
\end{itemize}

\section{The subalgebras of $\ct$}

\begin{figure}[!h]
\includegraphics[scale=0.85]{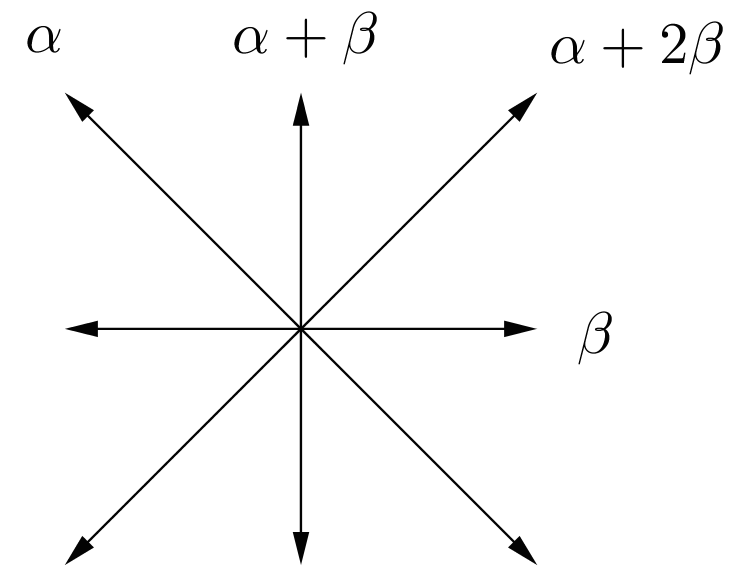}
\caption{The root system of $C_2$}\label{ct_roots_ref}
\end{figure}

The simple Lie algebra $C_2$ may be realized as the symplectic algebra $\mathfrak{sp}(4,\mathbb{C})$ 
of $4\times 4$ complex matrices $X$ satisfying $J X^{t}J=X$, where $J$ is the  matrix
\begin{equation}\label{eq:matrixJ}
  J = \begin{pmatrix}
    0&0&1&0\\ 
    0&0&0&1\\
    -1&0&0&0\\
    0&-1&0&0\\
    \end{pmatrix}.
\end{equation}
The root system associated with $C_2$ is depicted in Figure \ref{ct_roots_ref}.  

It has positive roots $\alpha, ~\beta, ~\alpha+\beta$, and $\alpha+2\beta$. Under the identification with $\mathfrak{sp}(4,\mathbb{C})$, $C_2$ has a Chevalley basis  
\begin{equation}
\svast\{
\begin{array}{lllllllllll}
 H_\alpha, &  H_\beta, \\
 X_\alpha, &  X_\beta, &X_{\alpha+\beta},  &X_{\alpha+2\beta},\\
 Y_\alpha, & Y_\beta, &  Y_{\alpha+\beta},  &Y_{\alpha+2\beta}
\end{array}\svast\},
\end{equation}
where
\begin{eqnarray}\label{eq:rootvectors}
X_{\alpha} = \begin{pmatrix}
0&0&0&0\\
0&0&0&1\\
0&0&0&0\\
0&0&0&0\\
\end{pmatrix}, 
&X_{\beta} = \begin{pmatrix}
0&1&0&0\\
0&0&0&0\\
0&0&0&0\\
0&0&-1&0\\
\end{pmatrix}, \\
X_{\alpha+\beta} = \begin{pmatrix} \label{eq:rootvectors2}
0&0&0&1\\
0&0&1&0\\
0&0&0&0\\
0&0&0&0\\
\end{pmatrix}, 
&X_{\alpha + 2\beta} = \begin{pmatrix}
0&0&1&0\\
0&0&0&0\\
0&0&0&0\\
0&0&0&0\\
\end{pmatrix},
\end{eqnarray}
and 
\begin{equation}
\begin{array}{lllllllllllll}
\displaystyle Y_\alpha &=& X_\alpha^t, & Y_\beta &=& X_\beta^t, \\ [3pt]
\displaystyle Y_{\alpha+\beta} &=& X_{\alpha+\beta}^t, & Y_{\alpha+2\beta} &=& X_{\alpha+2\beta}^t, \\ [3pt]
\displaystyle H_\alpha&=&[X_\alpha, Y_\alpha], & H_\beta&=&[X_\beta, Y_\beta].
\end{array}
\end{equation}
$H_\alpha$ and $H_\beta$ are the coroots; $X_\alpha, X_\beta, X_{\alpha+\beta}$, and $X_{\alpha+2\beta}$ are the positive root vectors; and 
$Y_\alpha, Y_\beta, Y_{\alpha+\beta}$, and $Y_{\alpha+2\beta}$ are the negative root vectors. The classification  becomes somewhat less cumbersome if we also make the definition:
\begin{equation}
T_{a,b}= (a+b)H_\alpha+aH_\beta.
\end{equation}

We highlight three distinguished subalgebras: A Cartan subalgebra $\mathfrak{h}=\langle H_\alpha, H_\beta \rangle$; a Borel subalgebra $\mathfrak{b}=\langle H_\alpha, H_\beta, X_\alpha, X_\beta, X_{\alpha+\beta}, X_{\alpha+2\beta} \rangle$; and the nilradical of the Borel subalgebra $\mathfrak{n}=\langle X_\alpha, X_\beta, X_{\alpha+\beta}, X_{\alpha+2\beta} \rangle$.


The Lie group corresponding to $C_2$, under its identification with $\mathfrak{sp}(4, \mathbb{C})$, is the symplectic group $\mathrm{Sp}(4, \mathbb{C})$. It is the Lie group of $4\times 4$ matrices $M$ with complex entries  such that $M J M^t=J$.

The classification and identification of the subalgebras of $C_2$ are described in tables below.  The classification is from \cite{dc2, c2}. The classification is organized by dimension, and further organized by structure. Tables \ref{ct1}-\ref{ct2} contain the classification and identification of the solvable subalgebras, and Table \ref{ct3} contains the classification and identification of the semisimple and Levi decomposable subalgebras.

In the following subsections,  we will  describe the techniques or strategies used in the classification of solvable subalgebras of $C_2$ and the
Levi decomposable subalgebras of $C_2$. We will only do this for the $C_2$ case, and we will omit virtually all proofs. We'll also 
discuss
identification of subalgebras with respect to the classification of solvable and Levi decomposable  Lie algebras in \cite{levi}.

\subsection{Classifying and identifying the solvable subalgebras of $C_2$}

We proceed by dimension to classify the solvable subalgebras of $\ct$. 
For 
each dimension, we  further subdivide the classification according to  structure.

In dimension $1$, we consider  three subcases depending on whether  the generator is semisimple, nilpotent, or has a non-trivial Jordan decomposition. 
If the generator is nilpotent, we may use the rich theory on nilpotent orbits that may be found in Collingwood and McGovern's classic text \cite{collingwood}.
Specifically, there  are precisely three nonzero nilpotent orbits of $\ct$ [\cite{collingwood}, Theorem 5.1.3] with representatives $X_\alpha, X_\beta$, and $X_\alpha +X_\beta$, giving us the three possible inequivalent subalgebras with nilpotent generator
\begin{equation}\label{odcad}
\begin{array}{llllllll}
\langle X_{\beta} \rangle, &\langle X_{\alpha} \rangle,&
\langle X_{\alpha}+X_{\beta} \rangle.
\end{array}
\end{equation}

In the case where the generator is semisimple, we first use the fact that 
every semisimple element $T$ of $C_2$ is conjugate to an element in $\mathfrak{h}$ [\cite{collingwood}, Corollary 2.2.2], so we may assume $T\in \mathfrak{h}$. We also use the result that two elements in $\mathfrak{t}$ are $\mathrm{Sp}(4,\mathbb{C})$-conjugate if and only if they are $W$-conjugate [\cite{collingwood}, Theorem 2.2.4], where $W$ is the Weyl group corresponding to $\mathfrak{h}$.  The Weyl group $W$ of $C_2$  has 
generators 
$s_\alpha$ and $s_\beta$ such that  $s_\alpha(T_{a,b})=T_{a,-b}$ and $s_\beta(T_{a,b})=T_{b,a}$. It follows that   a complete list of   one-dimensional subalgebras of $C_2$ with semisimple generators is
\begin{equation} \label{odcbd}
\begin{array}{llllllll}
\langle T_{1,b} \rangle \cong \langle T_{1,b^{-1}} \rangle,  b \neq 0,  \pm 1; &\langle T_{1,0} \rangle;&
\langle T_{1, 1}\rangle. 
\end{array}
\end{equation} 

We then rely primarily on linear algebra to consider the case where the generator has a non-trivial Jordan decomposition. That is, the generator
is the sum of a non-zero nilpotent element and a non-zero semisimple element 
that commutes with it.
This yields a complete list of  one-dimensional subalgebras of $C_2$ with generators  having a nontrivial Jordan decomposition:
\begin{equation}\label{opik6}
\begin{array}{llllllll}
\langle T_{1, 0}+X_{\alpha} \rangle; &\langle T_{1,1}+ X_{\beta} \rangle.
\end{array}
\end{equation}
Eqs. \eqref{odcad}-\eqref{opik6} give us a complete classification of one-dimensional subalgebras of $\ct$.

We proceed to higher dimensional solvable subalgebras, one dimension at a time, building from the previous dimension. Each dimension is further subdivided according to the nature of the Jordan decomposition of its generators, much as in  the one-dimensional case. When advantageous--specifically in dimensions three through five--we further differentiate elements that are regular and non-regular.  Following this general strategy, we get the classification of solvable subalgebras in Tables \ref{ct1}-\ref{ct2}.

To assist in the identification of  the solvable subalgebras (with respect to the partial classification of solvable Lie algebras described in \cite{levi}) found in the classification, we calculate certain data for each subalgebra. We often do this with the assistance of the 
computer algebra system Maple.  

The following is a list of useful data or characteristics that may be computed with Maple:  indecomposability; computation of the nilradical; computation of a basis for the Lie algebra of derivations; and the dimensions of the subalgebras in the derived series, upper central series, and lower central series. The aforementioned data may not uniquely identify the subalgebra in question, but it will narrow down the possibilities.  The generalized Casimir invariants and cohomology theory may also be employed to aid with identification. 

\subsection{Classifying and identifying the Levi decomposable subalgebras of $C_2$}

To determine the Levi decomposable  subalgebras of $\ct$, we first decompose $\ct$ with respect to the adjoint actions of each of its semisimple subalgebras. 
For instance, decomposing $\ct$ with respect to the simple subalgebra $\langle X_{\alpha+2\beta},$  $Y_{\alpha+2\beta}, T_{1,0} \rangle$ $\cong A_1$ yields

\begin{equation}\label{accf1}
\begin{array}{ccccccccccccc}
\ct &\cong_{A_1}& \langle X_{\alpha+2\beta},Y_{\alpha+2\beta}, T_{1,0}\rangle&\oplus& \langle X_\beta, Y_{\alpha+\beta} \rangle &\oplus&\\
&&\langle X_{\alpha+\beta}, Y_{\alpha+2\beta} \rangle&\oplus&\langle  X_\alpha \rangle &\oplus& \\
&&\langle Y_\alpha \rangle&\oplus&\langle  T_{0,1} \rangle  \\
&  \cong_{A_1} & V_{}(2) &\oplus&V_{}(1)&\oplus& \\
&&V_{}(1) &\oplus&V_{}(0)&\oplus& \\
&&V_{}(0) &\oplus&V_{}(0),
\end{array}
\end{equation}
where $V(n)$ is the $n+1$ dimensional, irreducible representation of $A_1$.  In the  decomposition  (the first instance of) $V(2)$ is isomorphic to $A_1$ as a subalgebra. The remaining representations in the decomposition, or combinations thereof, give us the potential radical components for Levi decomposable subalgebras of $\ct$, having $\langle X_{\alpha+2\beta}, Y_{\alpha+2\beta}, T_{1,0} \rangle$ as the Levi factor.  

We find that the following is a complete list of inequivalent subalgebras of $\ct$ having $\langle X_{\alpha+2\beta}, Y_{\alpha+2\beta}, T_{1,0} \rangle$ as the Levi factor:
\begin{equation}\label{tgh33}
\begin{array}{lllll}
\langle X_{\alpha+2\beta}, Y_{\alpha+2\beta}, T_{1,0} \rangle \oplus \langle X_\alpha \rangle,\\
 \langle X_{\alpha+2\beta}, Y_{\alpha+2\beta}, T_{1,0} \rangle \oplus \langle T_{0,1}\rangle,\\
 \langle X_{\alpha+2\beta}, Y_{\alpha+2\beta}, T_{1,0} \rangle \oplus \langle X_\alpha, T_{0,1} \rangle,\\
 \langle X_{\alpha+2\beta}, Y_{\alpha+2\beta}, T_{1,0} \rangle \inplus \langle X_\beta, Y_{\alpha+\beta}, Y_\alpha\rangle,\\
 \langle X_{\alpha+2\beta}, Y_{\alpha+2\beta}, T_{1,0} \rangle \inplus \langle X_\beta, Y_{\alpha+\beta}, Y_\alpha, T_{0,1}\rangle.
\end{array}
\end{equation}
We continue this procedure for each semisimple subalgebra of $\ct$ to obtain a classification of the Levi decomposable subalgebras of $\ct$.

To identify the Levi decomposable subalgebras with respect to the (partial) classification of Levi decomposable subalgebra described in \cite{levi},
we identify the radical, and the representation of the Levi factor on the radical. In this manner, we get the following identification of the subalgebras 
in Eq. \eqref{tgh33}:
\begin{equation}
\begin{array}{lllllll}
\langle X_{\alpha+2\beta}, Y_{\alpha+2\beta}, T_{1,0} \rangle \oplus \langle X_\alpha \rangle &\cong& A_1 \oplus \mathfrak{n}_{1,1},\\
 \langle X_{\alpha+2\beta}, Y_{\alpha+2\beta}, T_{1,0} \rangle \oplus \langle T_{0,1}\rangle&\cong& A_1 \oplus \mathfrak{n}_{1,1},\\
 \langle X_{\alpha+2\beta}, Y_{\alpha+2\beta}, T_{1,0} \rangle \oplus \langle X_\alpha, T_{0,1} \rangle&\cong& A_1 \oplus \mathfrak{s}_{2,1},\\
 \langle X_{\alpha+2\beta}, Y_{\alpha+2\beta}, T_{1,0} \rangle \inplus \langle X_\beta, Y_{\alpha+\beta}, Y_\alpha\rangle&\cong&  A_1 \inplus \mathfrak{n}_{3,1},\\
 \langle X_{\alpha+2\beta}, Y_{\alpha+2\beta}, T_{1,0} \rangle \inplus \langle X_\beta, Y_{\alpha+\beta}, Y_\alpha, T_{0,1}\rangle&\cong&  A_1 \inplus \mathfrak{s}_{4,8, A=1}.
\end{array}
\end{equation}
Precise descriptions of the Levi decomposable algebras listed above may be found in the appendices.


\begin{small}
\renewcommand{\arraystretch}{1.5}
\begin{table}[h!]
\begin{tabular}{|c|c|c|c|c|c|c|}
  \hline 
  Representative & Conditions & Equivalences &  Isomorphism Class\\ 
  \hline
 \multicolumn{4}{|c|}{Semisimple}  \\
  \hline
$\langle T_{a,b} \rangle$ & $a,b \ne 0$  & $\langle T_{a,b}    \rangle \sim \langle T_{a',b'} \rangle $   &$\mathfrak{n}_{1,1}$\\
&$a \ne \pm b$& iff $\{ a,b \}=\{ \lambda a',\pm \lambda b' \}$ & \\
  &&  for some $\lambda \in \mathbb{C}^*$&\\
  \hline
  $\langle T_{1,0} \rangle$ & & & $\mathfrak{n}_{1,1}$ \\
  \hline
$\langle T_{1,1} \rangle$ & & & $\mathfrak{n}_{1,1}$\\
\hline 
 \multicolumn{4}{|c|}{Nilpotent} \\
  \hline
  $\langle X_{\alpha} \rangle$ & &  & $\mathfrak{n}_{1,1}$\\  
      \hline
$\langle X_{\beta} \rangle$ & &  &$\mathfrak{n}_{1,1}$\\  
      \hline
$\langle X_{\alpha} +X_{\beta} \rangle$ & & &$\mathfrak{n}_{1,1}$ \\  
      \hline
  \multicolumn{4}{|c|}{  Non-Trivial Jordan Decomposition} \\
  \hline
  $\langle T_{1,0} + X_{\alpha}  \rangle $ & &
& $\mathfrak{n}_{1,1}$\\
\hline 
$\langle T_{1,1} + X_{\beta} \rangle$ &  &
  &  $\mathfrak{n}_{1,1}$ \\
\hline 
\end{tabular}
\vspace{1mm}
\caption{One-dimensional (solvable)  subalgebras of $C_2$}\label{onedimuu}\label{ct1}
\end{table}
\end{small}

\begin{small}
\renewcommand{\arraystretch}{1.5}
\begin{table}[h!]
\begin{tabular}{|c|c|c|c|c|c|}
  \hline 
  Representative & Conditions & Equivalences &  Isomorphism Class\\ 
  \hline
 \multicolumn{4}{|c|}{Semisimple}  \\
  \hline
   $\algt$ &&&2$\mathfrak{n}_{1,1}$ \\ 
  \hline
   \multicolumn{4}{|c|}{Containing both Semisimple and Nilpotent Elements}  \\
  \hline
     $\langle T_{3,1}, X_{\alpha} + X_{\beta} \rangle$ &&&$\mathfrak{s}_{2,1}$\\
\hline
$\langle T_{a,1} , X_{\alpha} \rangle$ & $  a \ne 0, \pm 1$  & $\langle T_{a,1} , X_{\alpha} \rangle \sim \langle T_{-a,1} , X_{\alpha} \rangle$ &$\mathfrak{s}_{2,1}$\\
\hline

$\langle T_{a,1} , X_{\beta} \rangle$ & $ a \ne 0, \pm 1$ &   $\langle T_{a,1} , X_{\beta} \rangle \sim \langle T_{a^{-1},1} , X_{\beta} \rangle$ &$\mathfrak{s}_{2,1}$\\
\hline
 $\langle  T_{1,0}, X_{\alpha} \rangle$ &&&2$\mathfrak{n}_{1,1}$\\ 
\hline
  $ \langle  T_{1,0}, X_{\beta} \rangle$ &&&$\mathfrak{s}_{2,1}$\\ 
\hline
 $\langle  T_{1,0}, X_{\alpha  +2 \beta} \rangle $ &&&$\mathfrak{s}_{2,1}$\\
\hline
 $\langle  T_{1,1}, X_{\alpha} \rangle$ &&&$\mathfrak{s}_{2,1}$\\
\hline 
 $\langle  T_{1,1}, X_{\beta} \rangle$ &&&2$\mathfrak{n}_{1,1}$\\
\hline 
 $\langle  T_{1,1}, X_{\alpha+ \beta} \rangle$ &&&$\mathfrak{s}_{2,1}$\\
\hline
  \multicolumn{4}{|c|}{Containing no Semisimple Elements, but not Nilpotent}  \\
  \hline
  $\langle T_{1,1}+X_{\beta}, X_{\alpha + 2 \beta} \rangle$ &&&$\mathfrak{s}_{2,1}$\\
\hline 
 $\langle T_{1,0}+X_{\alpha }, X_{\alpha + \beta } \rangle$&&&$\mathfrak{s}_{2,1}$\\
\hline 
 $\langle T_{1,0}+X_{\alpha }, X_{\alpha + 2\beta } \rangle$&&&$\mathfrak{s}_{2,1}$\\
\hline
  \multicolumn{4}{|c|}{Nilpotent}  \\
  \hline
  $\langle X_{\alpha}, X_{\alpha + \beta} \rangle $ &&& 2$\mathfrak{n}_{1,1}$\\
\hline 
 $\langle X_{\alpha} ,  X_{\alpha + 2 \beta}  \rangle $ &&&  2$\mathfrak{n}_{1,1}$\\
\hline 
 $\langle X_{\beta} + X_{\alpha}  , X_{\alpha+ 2 \beta} \rangle$ &&&  2$\mathfrak{n}_{1,1}$\\
\hline 
\end{tabular}
\vspace{1mm}
\caption{Two-dimensional (solvable) subalgebras of $C_2$}\label{twodimuu}
\end{table}
\end{small}


\begin{landscape}
\begin{scriptsize}
\renewcommand{\arraystretch}{1.4}
\begin{table}[h!]
\begin{tabular}{|c|c|c|c|c|c|}
  \hline 
  Representative & Conditions & Equivalences & Isomorphism Class\\ 
  \hline
 \multicolumn{4}{|c|}{Containing a Cartan Subalgebra}  \\
  \hline
$\langle \algt, X_{\alpha} \rangle$ &   & &  $\mathfrak{n}_{1,1}\oplus \mathfrak{s}_{2,1}$\\
  \hline
  $\langle \algt, X_{\beta} \rangle$ &   & & $\mathfrak{n}_{1,1}\oplus \mathfrak{s}_{2,1}$\\
  \hline
   \multicolumn{4}{|c|}{Containing a Regular Semisimple Element but not a Cartan Subalgebra}  \\
  \hline
   $\langle  T_{a,1} , X_{\alpha}, X_{\alpha+ \beta} \rangle$ &$ a \ne 0, \pm 1, -3$&&  $\mathfrak{s}_{3,1, A=\frac{1+2\alpha  + \sqrt{1+4\alpha}}{-2\alpha}}$\\
   &&&where $\alpha=-2(a+1)/(a+3)^{2}$\\
\hline 
 $\langle  T_{-3,1} , X_{\alpha}, X_{\alpha+ \beta} \rangle$ &&& $\mathfrak{s}_{3,1, A=-1}$\\
\hline 
 $\langle T_{a,1} , X_{\alpha}, X_{\alpha+ 2\beta} \rangle$  & $a \ne 0, \pm 1$ &  $\langle T_{a,1} , X_{\alpha}, X_{\alpha+ 2\beta} \rangle 
\sim  \langle T_{a^{-1},1} , X_{\alpha}, X_{\alpha+ 2\beta} \rangle$ &  $\mathfrak{s}_{3,1, A=\frac{1+2\alpha   + \sqrt{1+4\alpha}}{-2\alpha}}$\\
 &&& where $\alpha=- a/(a+1)^{2} $\\
\hline 
 $\langle T_{3,1} , X_{\alpha}+ X_{\beta},  X_{\alpha+ 2\beta} \rangle$ &&&   $\mathfrak{s}_{3,1, A=1/3}$ \\
  \hline
  \multicolumn{4}{|c|}{Containing a Non-regular Semisimple Element but not a Cartan Subalgebra}  \\      
      \hline
 $\langle  T_{1,0}, X_{\alpha}, X_{\alpha + \beta} \rangle$ &&& $\mathfrak{n}_{1,1}\oplus \mathfrak{s}_{2,1}$\\  
 \hline
 $\langle  T_{1,0}, X_{\alpha}, X_{\alpha + 2 \beta} \rangle$ &&&  $\mathfrak{n}_{1,1}\oplus \mathfrak{s}_{2,1}$  \\  
 \hline
 $\langle  T_{1,0}, X_{\alpha + \beta}, X_{\alpha + 2\beta} \rangle$ &&&    $\mathfrak{s}_{3,1, A=1/2}$\\  
\hline
 $\langle T_{1,-1} , X_{\alpha  + \beta},  X_{\alpha + 2 \beta} \rangle$ &&&    $\mathfrak{n}_{1,1}\oplus \mathfrak{s}_{2,1}$\\
 \hline
 $\langle T_{1,-1} , X_{\alpha},  X_{\alpha + 2\beta} \rangle$ &&&  $\mathfrak{s}_{3,1, A=-1}$\\
 \hline
 $\langle T_{1,-1} , X_{\beta} ,  X_{\alpha + 2\beta} \rangle$ &&&  $\mathfrak{s}_{3,1, A=1}$\\
 \hline
 $\langle T_{1,1} , X_{\alpha}  , X_{\alpha + 2 \beta} \rangle$ &&&$\mathfrak{s}_{3,1, A=1}$\\
\hline
  \multicolumn{4}{|c|}{Containing no Semisimple Elements but not Nilpotent}  \\      
      \hline
 $\langle T_{1,1} + X_{\beta}, X_{\alpha+\beta} , X_{\alpha+2\beta} \rangle$ &&&  $\mathfrak{s}_{3,2}$\\ 
\hline 
 $\langle T_{1,-1} + X_{\alpha+\beta} , X_{\alpha}, X_{\alpha+2\beta} \rangle$ &&&  $\mathfrak{s}_{3,1, A=-1}$ \\
\hline 
 $\langle T_{1,0} + X_{\alpha} , X_{\alpha + \beta}, X_{\alpha+2\beta} \rangle$ &&& $\mathfrak{s}_{3,1, A=1/2}$\\ 
\hline
  \multicolumn{4}{|c|}{Nilpotent }  \\
  \hline
  $\langle X_\alpha, X_{\alpha+\beta}, X_{\alpha+2\beta}\rangle$ & &  & $3\mathfrak{n}_{1,1}$\\  
      \hline
$\langle X_{\beta},  X_{\alpha + \beta}, X_{\alpha + 2 \beta} \rangle$ & &  &$\mathfrak{n}_{3,1}$\\  
      \hline
$\langle X_{\alpha} + X_{\beta},  X_{\alpha + \beta}, X_{\alpha + 2 \beta} \rangle$ & & & $\mathfrak{n}_{3,1}$\\  
      \hline
\end{tabular}
\vspace{1mm}
\caption{Three-dimensional solvable subalgebras of $C_2$}\label{threedimuu}
\end{table}
\end{scriptsize}
\end{landscape}


\begin{scriptsize}
\renewcommand{\arraystretch}{1.5}
\begin{table}[h!]
\begin{tabular}{|c|c|c|c|c|}
  \hline 
  Representative & Conditions & Equivalences & Isomorphism Class\\ 
  \hline
 \multicolumn{4}{|c|}{Containing a Cartan Subalgebra}  \\
  \hline
   $  \langle \algt, X_{\alpha}, X_{\alpha + \beta} \rangle $ &&& $\mathfrak{s}_{4,12}$\\ 
\hline
 $\langle \algt, X_{\alpha}, X_{\alpha + 2\beta} \rangle$ &&& $\mathfrak{s}_{4,12}$ \\ 
\hline
 \multicolumn{4}{|c|}{Containing a Regular Semisimple Element but not a Cartan Subalgebra}  \\
  \hline
   $    \langle T_{a,1}, X_\alpha, X_{\alpha+\beta}, X_{\alpha+2\beta} \rangle$ & $a \ne 0, \pm 1$ &  $\langle T_{a,1},X_\alpha, X_{\alpha+\beta}, X_{\alpha+2\beta} \rangle$ 
 &    $\mathfrak{s}_{4, 2}$, $\mathfrak{s}_{4,3}$, $\mathfrak{s}_{4,4}$,  \\
     &  & $\sim \langle T_{-a,1}, X_{\beta}, X_{\alpha + \beta}, X_{\alpha + 2 \beta} \rangle$ &  $\mathfrak{n}_{1,1}\oplus \mathfrak{s}_{3,1}$,
     or $\mathfrak{n}_{1,1}\oplus \mathfrak{s}_{3,2}$ \\
        && $\sim \langle T_{a^{-1},1}, X_\alpha, X_{\alpha+\beta}, X_{\alpha+2\beta} \rangle$&   
        (depending on $a$, c.f., \cite{c2})
        \\
  \hline
  $\langle T_{a,1}, X_{\beta}, X_{\alpha + \beta}, X_{\alpha + 2 \beta} \rangle$ & $a \ne 0, \pm 1$ & $\langle T_{a,1}, X_{\beta}, X_{\alpha + \beta}, X_{\alpha + 2 \beta} \rangle$ & $\mathfrak{s}_{4, 8, A= \frac{1+2\alpha   + \sqrt{1+4\alpha}}{-2\alpha}}$ \\
   &  & $\sim \langle T_{-a,1}, X_{\beta}, X_{\alpha + \beta}, X_{\alpha + 2 \beta} \rangle$ & where $\alpha=(1-a^{2})/(4a^{2})$ \\
  \hline
    $ \langle T_{3,1}, X_{\alpha} + X_{\beta} , X_{\alpha + \beta}, X_{\alpha + 2 \beta} \rangle$ &&  & $\mathfrak{s}_{4, 8,  A=1/2}$\\
\hline

  \multicolumn{4}{|c|}{Containing a Non-regular Semisimple Element but not a Cartan Subalgebra}  \\      
      \hline
 $\langle T_{0,1}, X_\alpha, X_{\alpha+\beta}, X_{\alpha+2\beta} \rangle$ &&& $\mathfrak{n}_{1,1}\oplus \mathfrak{s}_{3, 1, A=1/2}$ \\
\hline
 $\langle T_{0,1}, X_{\beta}, X_{\alpha + \beta}, X_{\alpha + 2 \beta} \rangle$ &&&   $\mathfrak{s}_{4,6}$ \\ 
\hline
 $\langle T_{1,0}, X_{\beta}, X_{\alpha + \beta}, X_{\alpha + 2 \beta} \rangle$ &&&  $\mathfrak{s}_{4,8, A=1}$\\
\hline
 $\langle T_{1,1}, X_\alpha, X_{\alpha+\beta}, X_{\alpha+2\beta} \rangle$ &&&   $\mathfrak{s}_{4, 3, A=B=1}$\\
\hline
 $\langle T_{1,-1}, X_\alpha, X_{\alpha+\beta}, X_{\alpha+2\beta} \rangle$ &&&  $\mathfrak{n}_{1,1}\oplus \mathfrak{s}_{3,1, A=-1}$  \\
\hline
 $\langle T_{1,1}, X_{\beta}, X_{\alpha + \beta}, X_{\alpha + 2 \beta} \rangle$ &&&   $\mathfrak{s}_{4, 11}$ \\
\hline
 \multicolumn{4}{|c|}{Containing no Semisimple Elements but not Nilpotent}  \\
  \hline
  $\langle T_{1,1} + X_{\beta}, X_{\alpha}, X_{\alpha+ \beta}, X_{\alpha+ 2\beta} \rangle$ &&&  $\mathfrak{s}_{4, 2}$\\
\hline
$\langle  T_{1,0} + X_{\alpha} , X_{\beta},  X_{\alpha+ \beta}, X_{\alpha+ 2\beta} \rangle$ &&& $\mathfrak{s}_{4, 10}$ \\
\hline
\multicolumn{4}{|c|}{Nilpotent Subalgebras}  \\
  \hline
  $ \algn$ &&&    $\mathfrak{n}_{4,1}$\\
\hline
 \end{tabular}
\vspace{1mm}
\caption{Four-dimensional solvable subalgebras of $C_2$}\label{fourdimuu}
\end{table}
\end{scriptsize}



\renewcommand{\arraystretch}{1.8}
\begin{table}[h!]
\begin{tabular}{|c|c|c|c|c|c|}
  \hline 
 Dimension & Representative & Conditions &  Isomorphism Class\\ 
  \hline
 \multicolumn{4}{|c|}{Containing a Cartan Subalgebra}  \\
  \hline
$5$  &  $  \langle \algt, X_\alpha, X_{\alpha+\beta}, X_{\alpha+2\beta} \rangle $ & & $\mathfrak{s}_{5,41, A=B=\frac{1}{2}}$\\ 
\hline
$5$ & $\langle \algt , X_{\beta}, X_{\alpha + \beta}, X_{\alpha + 2\beta} \rangle$ &&  $\mathfrak{s}_{5,44}$\\ 
\hline
 \multicolumn{4}{|c|}{Containing a Regular Semisimple Element but not a Cartan Subalgebra}  \\
  \hline
 $5$   &$    \langle T_{a,1}, \algn \rangle$ & $a \ne 0, \pm 1$    
 &  $\mathfrak{s}_{5, 35, A=\frac{2}{a-1}}$ \\
  \hline
   \multicolumn{4}{|c|}{Containing a Non-regular Semisimple Element but not a Cartan Subalgebra}  \\
  \hline
$5$  & $\langle T_{1,-1}, \algn \rangle$ &&    $\mathfrak{s}_{5, 35, A=-1}$\\
\hline
$5$& $\langle T_{1,1}, \algn \rangle$ &&     $\mathfrak{s}_{5, 37}$\\
\hline
$5$& $\langle T_{1,0}, \algn \rangle$ &&   $\mathfrak{s}_{5, 36}$  \\
\hline
$5$& $\langle T_{0,1}, \algn \rangle$ &&  $\mathfrak{s}_{5, 33}$ \\
\hline 
\hline
$6$ & $ \borel$ &&    $\mathfrak{s}_{6, 242}$\\
\hline
 \end{tabular}
\vspace{1mm}
\caption{Five- and six-dimensional solvable subalgebras of $C_2$}\label{ct2}
\end{table}

\begin{small}
\renewcommand{\arraystretch}{1.5}
\begin{table}[h!]
\begin{tabular}{|c|c|c|c|c|c|}
  \hline 
 Dimension & Representative &   Isomorphism Class\\ 
      \hline
   \multicolumn{3}{|c|}{Semisimple Subalgebras}  \\
  \hline
$3$  &  $  \langle X_{\alpha+2\beta}, Y_{\alpha+2\beta}, T_{1,0} \rangle $ &  $A_1$\\ 
\hline
$3$  &  $  \langle X_{\alpha+\beta}, Y_{\alpha+\beta}, T_{1,1} \rangle $ &  $A_1$\\ 
\hline
$3$  &  $  \langle X_{\alpha}+X_{\beta}, Y_{\alpha}+Y_{\beta}, T_{3,1} \rangle $ &  $A_1$\\ 
\hline
$6$  &  $ \langle X_{\alpha}, Y_{\alpha}, T_{0,1} \rangle \oplus \langle X_{\alpha+2\beta}, Y_{\alpha+2\beta}, T_{1,0} \rangle$ &  $A_1\times A_1$\\ 
    \hline
   \multicolumn{3}{|c|}{Levi Decomposable Subalgebras}  \\
  \hline
$4$  &  $ \langle X_{\alpha+\beta}, Y_{\alpha+\beta}, T_{1,1} \rangle \oplus \langle T_{1,-1}\rangle$ &  $A_1\oplus \mathfrak{n}_{1,1}$\\ 
\hline
$4$  &  $  \langle X_{\alpha+2\beta}, Y_{\alpha+2\beta}, T_{1,0} \rangle \oplus  \langle X_{\alpha}\rangle$&  $A_1\oplus \mathfrak{n}_{1,1}$\\
\hline
$4$  &  $  \langle X_{\alpha+2\beta}, Y_{\alpha+2\beta}, T_{1,0} \rangle \oplus  \langle T_{0,1}\rangle$& $A_1\oplus \mathfrak{n}_{1,1}$\\ 
\hline
$5$  &  $  \langle X_{\alpha+2\beta}, Y_{\alpha+2\beta}, T_{1,0} \rangle \oplus  \langle X_{\alpha}, T_{0,1}\rangle$ &  $A_1\oplus  \mathfrak{s}_{2,1}$\\ 
\hline
$6$  & $ \langle X_{\alpha+\beta}, Y_{\alpha+\beta}, T_{1,1} \rangle \inplus \langle X_{\alpha}, Y_{\beta}, Y_{\alpha+2\beta}\rangle$  &$A_1\inplus 3\mathfrak{n}_{1,1}$ \\ 
\hline
$6$  & $  \langle X_{\alpha+2\beta}, Y_{\alpha+2\beta}, T_{1,0} \rangle \inplus  \langle X_{\beta}, Y_{\alpha+\beta}, Y_{\alpha}\rangle$ &  $A_1\inplus \mathfrak{n}_{3,1}$\\ 
\hline
$7$  & $ \langle X_{\alpha+\beta}, Y_{\alpha+\beta}, T_{1,1} \rangle \inplus \langle X_{\alpha}, Y_{\beta}, Y_{\alpha+2\beta}, T_{1,-1}\rangle$   & $A_1\inplus  \mathfrak{s}_{4,3, A=B=1}$\\ 
\hline
$7$  & $\langle X_{\alpha+2\beta}, Y_{\alpha+2\beta}, T_{1,0} \rangle \inplus \langle X_{\beta}, Y_{\alpha+\beta}, Y_{\alpha}, T_{0,1}\rangle$   & $A_1\inplus  \mathfrak{s}_{4,8, A=1}$ \\ 
\hline
 \end{tabular}
\vspace{1mm}
\caption{Semisimple and Levi decomposable subalgebras of $C_2$}\label{ct3}
\end{table}
\end{small}

\clearpage

\section{The subalgebras of $\at$}

\begin{figure}[!h]
\includegraphics[scale=0.85]{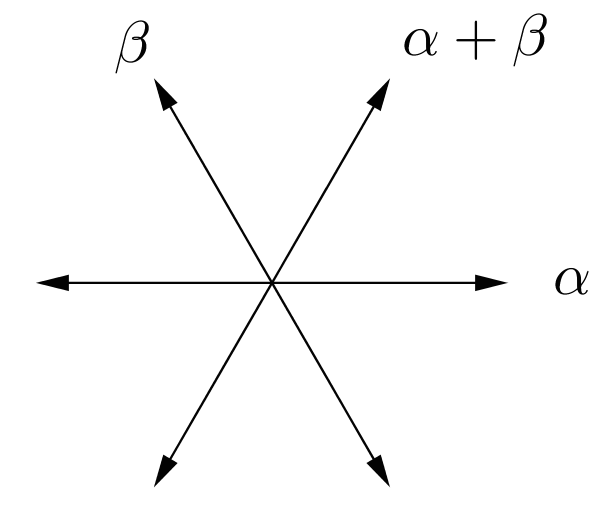}
\caption{The root system of $A_2$}\label{at_roots_ref}
\end{figure}

The simple Lie algebra $A_2$ may be realized as the special linear algebra $\mathfrak{sl}(3,\mathbb{C})$ of traceless
$3\times 3$ matrices with complex entries. The root system associated with $A_2$ is depicted in Figure \ref{at_roots_ref}.  

It has positive roots $\alpha, ~\beta$,  and $\alpha+\beta$. Under the identification with $\mathfrak{sl}(3,\mathbb{C})$, $A_2$ has a Chevalley basis  
\begin{equation}
\vast\{
\arraycolsep=1.7pt\def\arraystretch{1.2}
\begin{array}{llllllll}
\displaystyle H_\alpha, &  H_\beta, \\ 
 \displaystyle X_\alpha, &  X_\beta, &X_{\alpha+\beta}, \\ 
 \displaystyle  Y_\alpha, & Y_\beta, &  Y_{\alpha+\beta}
\end{array}\vast\},
\end{equation}
where 
\begin{equation}
\begin{array}{lllllllllllll}
H_\alpha&=& \left(
\begin{array}{rrr}
1&0&0 \\
0&-1&0 \\
0&0&0 
\end{array}
\right), &H_\beta&=& \left(
\begin{array}{rrr}
0&0&0 \\
0&1&0 \\
0&0&-1 
\end{array}
\right),\\ \\
X_\alpha&=& \left(
\begin{array}{rrr}
0&1&0 \\
0&0&0 \\
0&0&0 
\end{array}
\right), &Y_\alpha&=& \left(
\begin{array}{rrr}
0&0&0 \\
1&0&0 \\
0&0&0 
\end{array}
\right), 
\\ \\
X_\beta&=& \left(
\begin{array}{rrr}
0&0&0 \\
0&0&1 \\
0&0&0 
\end{array}
\right), &Y_\beta&=& \left(
\begin{array}{rrr}
0&0&0 \\
0&0&0 \\
0&1&0 
\end{array}
\right),\\ \\
X_{\alpha+\beta}&=& \left(
\begin{array}{rrr}
0&0&-1 \\
0&0&0 \\
0&0&0 
\end{array}
\right), &Y_{\alpha+\beta}&=& \left(
\begin{array}{rrr}
0&0&0 \\
0&0&0 \\
-1&0&0 
\end{array}
\right).
\end{array}
\end{equation}
$H_\alpha$ and $H_\beta$ are the coroots; $X_\alpha, X_\beta$, and $X_{\alpha+\beta}$ are the positive root vectors; and 
$Y_\alpha, Y_\beta$, and $Y_{\alpha+\beta}$ are the negative root vectors.

We highlight three distinguished subalgebras: A Cartan subalgebra $\mathfrak{h}=\langle H_\alpha, H_\beta \rangle$; a Borel subalgebra $\mathfrak{b}=\langle H_\alpha, H_\beta, X_\alpha, X_\beta, X_{\alpha+\beta} \rangle$; and the nilradical of the Borel subalgebra $\mathfrak{n}=\langle X_\alpha, X_\beta, X_{\alpha+\beta} \rangle$.


The Lie group corresponding to $A_2$, under its identification with $\mathfrak{sl}(3, \mathbb{C})$, is the special linear group $\mathrm{SL}(3, \mathbb{C})$. It is the Lie group of $3\times 3$ matrices with complex entries
and determinant $1$.

The classification and identification of the subalgebras of $A_2$ are described in the following tables.  The classification is from \cite{a2}, and the identification was 
established as part of the present article. The classification is organized by dimension. Tables \ref{a21}-\ref{a2fs} contain the classification and identification of the solvable subalgebras, and Table \ref{a2l} contains the classification and identification of the semisimple and Levi decomposable subalgebras.

Note that in Table \ref{a2fs} below, we addressed one oversight from the previous version of this article \cite{rank2v5}.
Specifically, the previous version contained the following two four-dimensional, solvable subalgebras: 
$\big\langle X_\alpha, X_\beta$, $X_{\alpha+\beta}$, $\frac{1}{2}H_\alpha+H_\beta  \big\rangle$ and $\langle X_\beta, Y_\alpha, Y_{\alpha+\beta}$, $2H_\alpha+H_\beta \rangle$.  And, these subalgebras
are conjugate in $\mathrm{SL}_3(\mathbb{C})$. In particular, $G\big\langle X_\alpha, X_\beta, X_{\alpha+\beta},$ $\frac{1}{2}H_\alpha+H_\beta \big\rangle G^{-1}=\langle X_\beta, Y_\alpha, Y_{\alpha+\beta}, 2H_\alpha+H_\beta \rangle$, where
\begin{equation}
G=\left(
\begin{array}{rrr}
0&0&1 \\
-1&0&0 \\
0&-1&0 
\end{array}
\right)\in \mathrm{SL}(3,\mathbb{C}).
\end{equation}
 The equivalence of these subalgebras seems to also have been noted in \cite{chap}. In the table below, we exclude the later subalgebra.  For further comments on the classification in
 \cite{chap} see \cite{a2b} (cf. \cite{a2}).

\begin{table}[h!]
\scalebox{0.9}{
\begin{tabular}{|c|c|c|c|c|c|}
  \hline 
 Dimension & Representative &   Isomorphism Class\\ 
  \hline
$1$  & $\langle X_\alpha+X_\beta\rangle$  &  $\mathfrak{n}_{1,1}$\\ 
\hline
$1$  & $\langle X_\alpha\rangle$  &  $\mathfrak{n}_{1,1}$\\ 
\hline
$1$  & $\langle X_\alpha+H_\alpha+2H_\beta\rangle$  &  $\mathfrak{n}_{1,1}$\\ 
\hline
$1$  & $\langle H_\alpha+aH_\beta\rangle$  &  $\mathfrak{n}_{1,1}$\\ 
  & $\langle H_\alpha+aH_\beta\rangle \sim \langle H_\alpha+bH_\beta\rangle$  & \\ 
    & iff $b=a, \frac{1}{a}, 1-a, \frac{1}{1-a}, \frac{a}{a-1},$ or $\frac{a-1}{a}$  & \\ 
\hline
$2$  &  $\langle X_\alpha+X_\beta, X_{\alpha+\beta}\rangle$ & $2\mathfrak{n}_{1,1}$\\ 
\hline
$2$  &  $\langle X_\alpha, H_\alpha+2H_\beta\rangle$ &  $2\mathfrak{n}_{1,1}$\\ 
\hline
$2$  &  $\langle X_\alpha, X_{\alpha+\beta}\rangle$ &  $2\mathfrak{n}_{1,1}$\\ 
\hline
$2$  &  $\langle X_\alpha, Y_\beta\rangle$ &  $2\mathfrak{n}_{1,1}$\\ 
\hline
$2$  &  $\langle H_\alpha, H_\beta\rangle$ &  $2\mathfrak{n}_{1,1}$\\ 
\hline
$2$  &  $\langle X_\alpha+X_\beta, H_\alpha+H_\beta\rangle$ & $\mathfrak{s}_{2,1}$\\ 
\hline
$2$  &  $\langle X_\alpha, -H_\alpha+H_\beta+3X_{\alpha+\beta}\rangle$ & $\mathfrak{s}_{2,1}$\\ 
\hline
$2$  &  $\langle X_\alpha,  -2H_\alpha-H_\beta+3Y_\beta\rangle$ & $\mathfrak{s}_{2,1}$\\ 
\hline
$2$  &  $\langle X_\alpha, a H_\alpha+(2a+1)H_\beta\rangle$ & $\mathfrak{s}_{2,1}$\\ 
  &  $\langle X_\alpha, a H_\alpha+(2a+1)H_\beta\rangle \sim \langle X_\alpha, b H_\alpha+(2b+1)H_\beta\rangle$ & \\ 
  &iff $a=b$&\\
\hline
 \end{tabular}}
\vspace{1mm}
\caption{One- and two-dimensional (solvable) subalgebras of $A_2$}\label{a21}
\end{table}

\begin{table}[h!]
\scalebox{0.85}{
\begin{tabular}{|c|c|c|c|c|c|c|}
  \hline 
 Dimension & Representative & Isomorphism Class & Conditions\\ 
  \hline
$3$  &  $\langle X_\alpha, X_{\alpha+\beta}, 2H_\alpha+H_\beta\rangle$ &$\mathfrak{s}_{3,1, A=1}$&\\ 
\hline
$3$  &  $\langle X_\alpha, Y_\beta, H_\alpha-H_\beta\rangle$ &$\mathfrak{s}_{3,1, A=1}$& \\ 
\hline
$3$  &  $\langle X_\alpha, X_{\alpha+\beta}, 2H_\alpha+H_\beta+X_\beta\rangle$ &  $\mathfrak{s}_{3,2}$&\\ 
\hline
$3$  &  $\langle Y_\alpha, Y_{\alpha+\beta}, 2H_\alpha+H_\beta+X_\beta\rangle$ & $\mathfrak{s}_{3,2}$&\\ 
\hline
$3$  &  $\langle X_\alpha+X_\beta, X_{\alpha+\beta}, H_\alpha+H_\beta\rangle$ &$\mathfrak{s}_{3,1, A=1/2}$&\\ 
\hline
$3$  &  $\langle X_\alpha, H_\alpha, H_\beta\rangle$ &$\mathfrak{n}_{1,1}\oplus \mathfrak{s}_{2,1}$& \\ 
\hline
$3$  &  $\mathfrak{u}^a=\langle X_\alpha, X_{\alpha+\beta}, (a-1)H_\alpha+aH_\beta\rangle$ &  $\mathfrak{n}_{1,1}\oplus \mathfrak{s}_{2,1}$&  if $a=\frac{1}{2}, 2$\\ 
&$a\neq \pm1$&   $\mathfrak{s}_{3,1, \frac{1+2\alpha+\sqrt{1+4\alpha}}{-2\alpha}}$,& otherwise\\
& $\mathfrak{u}^a \sim \mathfrak{u}^b$ iff $a=b$ or $ab=1$&$\alpha= -\frac{(2a-1)(a-2)}{9(a-1)^2}$ &\\
\hline
$3$  &  $\mathfrak{v}^a=\langle X_\alpha, Y_\beta, H_\alpha+aH_\beta\rangle$ &$\mathfrak{n}_{1,1}\oplus \mathfrak{s}_{2,1}$& if $a=\frac{1}{2}, 2$\\ 
&$a\neq \pm1$& $\mathfrak{s}_{3,1, \frac{1+2\alpha+\sqrt{1+4\alpha}}{-2\alpha}}$,& otherwise\\
&$\mathfrak{v}^a \sim \mathfrak{v}^b$ iff $a=b$ or $ab=1$& $\alpha= -\frac{(2a-1)(a-2)}{9(a-1)^2}$ &\\
\hline
$3$  &  $\langle X_\alpha, X_{\alpha+\beta}, H_\beta\rangle$ & $\mathfrak{s}_{3,1, A=-1}$&\\ 
\hline
$3$  &  $\langle X_\alpha, Y_\beta, H_\alpha+H_\beta\rangle$ & $\mathfrak{s}_{3,1, A=-1}$&\\ 
\hline
$3$  &  $\langle X_\alpha, X_\beta, X_{\alpha+\beta}\rangle$ &$\mathfrak{n}_{3,1}$& \\ 
\hline
$4$  & $\langle X_\alpha, X_{\alpha+\beta}, H_\alpha, H_\beta\rangle$  & $\mathfrak{s}_{4,12}$&\\ 
\hline
$4$  & $\langle X_\alpha, Y_\beta, H_\alpha, H_\beta\rangle$  & $\mathfrak{s}_{4,12}$& \\ 
\hline
$4$  & $\langle X_\alpha, X_\beta, X_{\alpha+\beta}, H_\alpha+H_\beta\rangle$  & $\mathfrak{s}_{4,8, A=1}$&\\ 
\hline
$4$  & $\langle X_\alpha, X_\beta, X_{\alpha+\beta},a H_\alpha+H_\beta\rangle$  & $\mathfrak{s}_{4, 11}$& if $a=\frac{1}{2}, 2$\\ 
&$a\neq \pm 1$& & \\
&$\langle X_\alpha, X_\beta, X_{\alpha+\beta},a H_\alpha+H_\beta\rangle \sim \langle X_\alpha, X_\beta, X_{\alpha+\beta}, b H_\alpha+H_\beta\rangle$& $\mathfrak{s}_{4, 8, \frac{1+2\alpha+\sqrt{1+4\alpha}}{-2\alpha}}$,&otherwise\\
&iff $a=b$& $\alpha=\frac{(2a-1)(a-2)}{(a+1)^2}$ &\\
\hline
$4$  & $\langle X_\alpha, X_\beta, X_{\alpha+\beta}, H_\alpha\rangle$  &   $\mathfrak{s}_{4, 8, A=-2}$&\\ 
\hline
$4$  & $\langle X_\alpha, X_\beta, X_{\alpha+\beta}, H_\alpha- H_\beta\rangle$  & $\mathfrak{s}_{4,6}$& \\ 
\hline
$5$  &  $\langle X_\alpha, X_\beta, X_{\alpha+\beta}, H_\alpha, H_\beta\rangle$&  $\mathfrak{s}_{5, 44}$&\\ 
\hline
 \end{tabular}}
\vspace{1mm}
\caption{Three-, four- and five-dimensional solvable subalgebras of $A_2$}\label{a2fs}
\end{table}


\renewcommand{\arraystretch}{1.5}
\begin{table}[h!]
\begin{tabular}{|c|c|c|c|c|c|}
  \hline 
 Dimension & Representative &   Isomorphism Class\\ 
    \hline
   \multicolumn{3}{|c|}{Semisimple Subalgebras}  \\
   \hline
$3$  &  $  \langle X_{\alpha+\beta}, Y_{\alpha+\beta}, H_\alpha+H_\beta \rangle$ &  $A_1$\\ 
\hline
$3$  &  $  \langle X_\alpha+X_\beta, 2Y_\alpha+2Y_\beta, 2H_\alpha+2H_\beta \rangle$  &  $A_1$\\ 
\hline
   \multicolumn{3}{|c|}{Levi Decomposable Subalgebras}  \\
  \hline
$4$  &  $  \langle X_{\alpha+\beta}, Y_{\alpha+\beta}, H_\alpha+H_\beta \rangle\oplus \langle H_\alpha-H_\beta\rangle$ & $A_1 \oplus \mathfrak{n}_{1,1}$\\ 
\hline
$5$  &  $  \langle X_{\alpha+\beta}, Y_{\alpha+\beta}, H_\alpha+H_\beta  \rangle\inplus \langle X_\alpha, Y_\beta\rangle$ &  $A_1 \inplus 2\mathfrak{n}_{1,1}$\\ 
\hline
$5$  &  $  \langle X_{\alpha+\beta}, Y_{\alpha+\beta}, H_\alpha+H_\beta \rangle\inplus \langle X_\beta, Y_\alpha\rangle$ & $A_1 \inplus 2\mathfrak{n}_{1,1}$ \\ 
\hline
$6$  &  $  \langle X_{\alpha+\beta}, Y_{\alpha+\beta}, H_\alpha+H_\beta  \rangle\inplus \langle X_\alpha, Y_\beta, H_\alpha-H_\beta\rangle$ & $A_1 \inplus \mathfrak{s}_{3,1, A=1}$ \\ 
\hline
$6$  &  $  \langle X_{\alpha+\beta}, Y_{\alpha+\beta}, H_\alpha+H_\beta  \rangle\inplus \langle X_\beta, Y_\alpha, H_\alpha-H_\beta\rangle$   &  $A_1 \inplus \mathfrak{s}_{3,1, A=1}$\\ 
\hline
 \end{tabular}
\vspace{1mm}
\caption{Semisimple and Levi decomposable subalgebras of $A_2$}\label{a2l}
\end{table}

\clearpage

\section{The subalgebras of $\aoao$}

\begin{figure}[!h]
\includegraphics[scale=0.85]{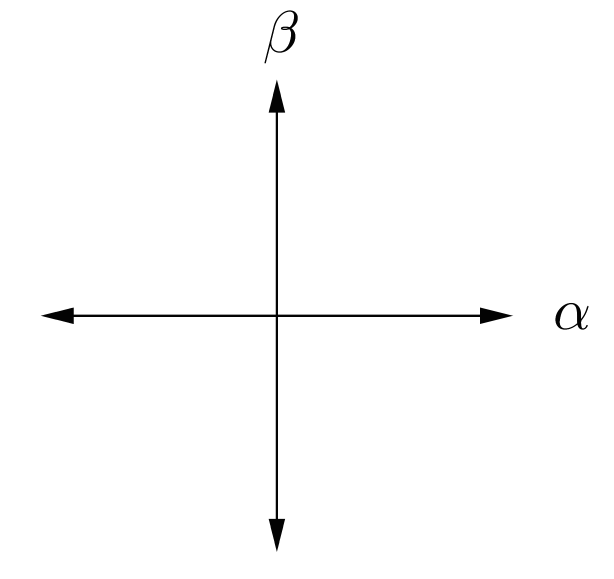}
\caption{The root system of $A_1\times A_1$}\label{aoao_rootss}
\end{figure}

The semisimple, and not simple, Lie algebra $A_1\times A_1$ may be realized as the special orthogonal algebra $\mathfrak{so}(4,\mathbb{C})$ of 
$4\times 4$ complex matrices satisfying $X^t=-X$. It is isomorphic to the
semisimple Lie algebra $\s \oplus \s$, where $\s$ is the special linear algebra of traceless $2\times 2$ matrices with complex entries.
The root system associated with $A_1\times A_1$ is depicted in Figure \ref{aoao_rootss}.  

It has positive roots $\alpha$ and $\beta$. Under the identification with $\s \oplus \s$, $A_1\times A_1$ has a Chevalley basis  
\begin{equation}
\begin{array}{llllllllllllll}
\{ H_\alpha, ~H_\beta, ~
 X_\alpha, ~X_\beta,  ~ Y_\alpha,  ~Y_\beta\}, 
 \end{array}
\end{equation}
where
\begin{equation}
\begin{array}{llllllllll}
H_\alpha&=& \left(
\begin{array}{rrrr}
1&0&0&0 \\
0&-1&0 &0\\
0&0&0 &0\\
0&0&0&0
\end{array}
\right), &H_\beta&=& \left(
\begin{array}{rrrr}
0&0&0& 0\\
0&0&0 &0\\
0&0&1 &0\\
0&0&0&-1
\end{array}
\right),  \\\\
X_\alpha&=& \left(
\begin{array}{rrrr}
0&1&0&0 \\
0&0&0 &0\\
0&0&0 &0\\
0&0&0&0
\end{array}
\right), &Y_\alpha&=& \left(
\begin{array}{rrrr}
0&0&0& 0\\
1&0&0 &0\\
0&0&0 &0\\
0&0&0&0
\end{array}
\right), \\\\
X_\beta&=& \left(
\begin{array}{rrrr}
0&0&0&0 \\
0&0&0 &0\\
0&0&0 &1\\
0&0&0&0
\end{array}
\right), &Y_\beta&=& \left(
\begin{array}{rrrr}
0&0&0& 0\\
0&0&0 &0\\
0&0&0 &0\\
0&0&1&0
\end{array}
\right).
\end{array}
\end{equation}
$H_\alpha$ and $H_\beta$ are the coroots; $X_\alpha$ and $X_{\beta} $ are the positive root vectors; and 
$Y_\alpha$ and $Y_{\beta}$ are the negative root vectors.

We highlight three distinguished subalgebras: A Cartan subalgebra $\mathfrak{h}=\langle H_\alpha, H_\beta \rangle$; a Borel subalgebra $\mathfrak{b}=\langle H_\alpha, H_\beta, X_\alpha, X_\beta \rangle$; and the nilradical of the Borel subalgebra $\mathfrak{n}=\langle X_\alpha, X_\beta \rangle$.

The Lie group corresponding to $A_1\times A_1$, under its identification with $\mathfrak{so}(4,\mathbb{C})$, is the special orthogonal group $\mathrm{SO}(4, \mathbb{C})$.  It is the Lie group of  $4 \times 4$ complex matrices $M$ with determinant $1$ satisfying $M^t M=I$.

The classification and identification of the subalgebras of $A_1\times A_1$ are described in the following tables.  The classification is from \cite{drsof}, and the identification was 
established 
as part of the present article. The classification is organized by dimension. Table \ref{dfg1} contains the classification and identification of the solvable subalgebras, and Table \ref{dfg2} contains the classification and identification of the semisimple and Levi decomposable subalgebras.

\begin{footnotesize}
\renewcommand{\arraystretch}{1.5}
\begin{table}[h!]
\begin{tabular}{|c|c|c|c|c|c|}
  \hline 
 Dimension & Representative  &  Isomorphism Class\\ 
  \hline
$1$  &  $  \langle X_\alpha \rangle$ &  $\mathfrak{n}_{1,1}$\\ 
\hline
$1$  &  $  \langle X_\beta \rangle$ &  $\mathfrak{n}_{1,1}$\\ 
\hline
$1$  &  $  \langle H_\alpha \rangle$ &  $\mathfrak{n}_{1,1}$\\ 
\hline
$1$  &  $  \langle H_\beta \rangle$ &  $\mathfrak{n}_{1,1}$\\ 
\hline
$1$  &  $  \langle X_\alpha+X_\beta \rangle$ &  $\mathfrak{n}_{1,1}$\\ 
\hline
$1$  &  $  \langle X_\alpha +H_\beta\rangle$ &  $\mathfrak{n}_{1,1}$\\ 
\hline
$1$  &  $  \langle H_\alpha +X_\beta\rangle$ &  $\mathfrak{n}_{1,1}$\\ 
\hline
$1$  &  $  \langle H_\alpha+aH_\beta \rangle$, ~~$a\in \mathbb{C}^*$&  $\mathfrak{n}_{1,1}$\\ 
  &  $ \langle H_\alpha+aH_\beta \rangle \sim  \langle H_\alpha+bH_\beta \rangle$  iff $a=\pm b$&  \\ 
\hline
$2$  &  $  \langle X_\alpha, X_\beta\rangle$ &  $2\mathfrak{n}_{1,1}$\\ 
\hline
$2$  & $  \langle X_\alpha, H_\beta\rangle$ &  $2\mathfrak{n}_{1,1}$\\ 
\hline
$2$  &  $  \langle H_\alpha, X_\beta\rangle$ &  $2\mathfrak{n}_{1,1}$\\ 
\hline
$2$  &  $  \langle H_\alpha, H_\beta\rangle$ &  $2\mathfrak{n}_{1,1}$\\ 
\hline
$2$  &  $  \langle X_\alpha+X_\beta, H_\alpha+H_\beta \rangle$ &  $\mathfrak{s}_{2,1}$\\ 
\hline
$2$  &  $  \langle X_\alpha, H_\alpha+aH_\beta \rangle$, ~~$a\in \mathbb{C}$ &  $\mathfrak{s}_{2,1}$\\ 
 &  $\langle X_\alpha, H_\alpha+aH_\beta \rangle\sim \langle X_\alpha, H_\alpha+bH_\beta \rangle$~ iff ~$a=\pm b$& \\ 
\hline
$2$  &  $  \langle X_\alpha, H_\alpha+X_\beta \rangle$ &  $\mathfrak{s}_{2,1}$\\ 
\hline
$2$  & $  \langle X_\beta, H_\beta+aH_\alpha \rangle$, ~~$a\in \mathbb{C}$ &  $\mathfrak{s}_{2,1}$\\ 
 &  $\langle X_\beta, H_\beta+aH_\alpha \rangle\sim \langle X_\beta, H_\beta+bH_\alpha \rangle$ ~iff ~$a=\pm b$& \\ 
\hline
$2$  &  $  \langle X_\beta, H_\beta+X_\alpha \rangle$ &  $\mathfrak{s}_{2,1}$\\ 
  \hline
$3$  &  $  \langle X_\alpha, X_\beta, H_\alpha+H_\beta \rangle$ &  $\mathfrak{s}_{3,1, A=1}$  \\ 
\hline
$3$  &  $  \langle X_\alpha, X_\beta, (1+\sqrt{1+4a})H_\alpha+(1-\sqrt{1+4a})H_\beta \rangle$ &    $\mathfrak{s}_{3,1, \frac{1+2a+\sqrt{1+4a}}{-2a}}$\\ 
  &  $a\in \mathbb{C}\setminus\{-\frac{1}{4},0\}$ &  \\ 
\hline
$3$  &  $  \langle X_\alpha, X_\beta, (1-\sqrt{1+4a})H_\alpha+(1+\sqrt{1+4a})H_\beta \rangle$ &   $\mathfrak{s}_{3,1, \frac{1+2a+\sqrt{1+4a}}{-2a}}$ \\ 
 &  $a\in \mathbb{C} \setminus \{-\frac{1}{4},0\}$ &  \\ 
\hline
$3$  &  $  \langle X_\beta, H_\alpha, H_\beta \rangle$ &  $\mathfrak{n}_{1,1}\oplus \mathfrak{s}_{2,1}$\\ 
\hline
$3$  &  $  \langle X_\alpha, X_\beta, H_\beta \rangle$ &  $\mathfrak{n}_{1,1}\oplus \mathfrak{s}_{2,1}$\\ 
\hline
$3$  &  $  \langle X_\alpha, H_\alpha, H_\beta \rangle$ &  $\mathfrak{n}_{1,1}\oplus \mathfrak{s}_{2,1}$\\ 
\hline
$3$  &  $  \langle X_\alpha, X_\beta, H_\alpha \rangle$ &  $\mathfrak{n}_{1,1}\oplus \mathfrak{s}_{2,1}$\\ 
\hline
$3$  &  $  \langle X_\alpha, X_\beta, H_\alpha-H_\beta \rangle$ &  $\mathfrak{s}_{3,1, A=-1}$\\ 
\hline
$4$  &  $  \langle X_\alpha, X_\beta, H_\alpha, H_\beta \rangle$ &  $\mathfrak{s}_{4,12}$\\ 
\hline
 \end{tabular}
\vspace{1mm}
\caption{Solvable subalgebras of $A_1\times A_1$}\label{dfg1}
\end{table}
\end{footnotesize}

\renewcommand{\arraystretch}{1.5}
\begin{table}[h!]
\begin{tabular}{|c|c|c|c|c|c|}
  \hline 
 Dimension & Representative &   Isomorphism Class\\ 
    \hline
   \multicolumn{3}{|c|}{Semisimple Subalgebras}  \\
  \hline
$3$  &  $  \langle X_\alpha, Y_\alpha, H_\alpha \rangle$ &  $A_1$\\ 
\hline
$3$  &  $  \langle X_\beta, Y_\beta, H_\beta \rangle$  &  $A_1$\\ 
\hline
$3$  &  $  \langle X_\alpha+X_\beta, Y_\alpha+Y_\beta,H_\alpha+H_\beta \rangle$  &  $A_1$\\ 
\hline
   \multicolumn{3}{|c|}{Levi Decomposable Subalgebras}  \\
  \hline
$4$  &  $  \langle X_\alpha, Y_\alpha, H_\alpha \rangle\oplus \langle H_\beta\rangle$ &  $A_1\oplus \mathfrak{n}_{1,1}$\\ 
\hline
$4$  &  $  \langle X_\alpha, Y_\alpha, H_\alpha \rangle\oplus \langle X_\beta\rangle$ &  $A_1\oplus \mathfrak{n}_{1,1}$\\ 
\hline
$4$  &  $  \langle X_\beta, Y_\beta, H_\beta \rangle\oplus \langle H_\alpha\rangle$ &  $A_1\oplus \mathfrak{n}_{1,1}$\\ 
\hline
$4$  &  $  \langle X_\beta, Y_\beta, H_\beta \rangle\oplus \langle X_\alpha\rangle$ &  $A_1\oplus \mathfrak{n}_{1,1}$\\ 
\hline
$5$  &  $  \langle X_\alpha, Y_\alpha, H_\alpha \rangle\oplus \langle X_\beta, H_\beta\rangle$   &  $A_1\oplus \mathfrak{s}_{2,1}$\\ 
\hline
$5$  &   $  \langle X_\beta, Y_\beta, H_\beta \rangle\oplus \langle X_\alpha, H_\alpha\rangle$  &  $A_1\oplus \mathfrak{s}_{2,1}$\\ 
\hline
 \end{tabular}
\vspace{1mm}
\caption{Semisimple and Levi decomposable subalgebras of $A_1\times A_1$}\label{dfg2}
\end{table}

\clearpage

\section{The subalgebras of $G_2$}

\begin{figure}[!h]
\includegraphics[scale=0.85]{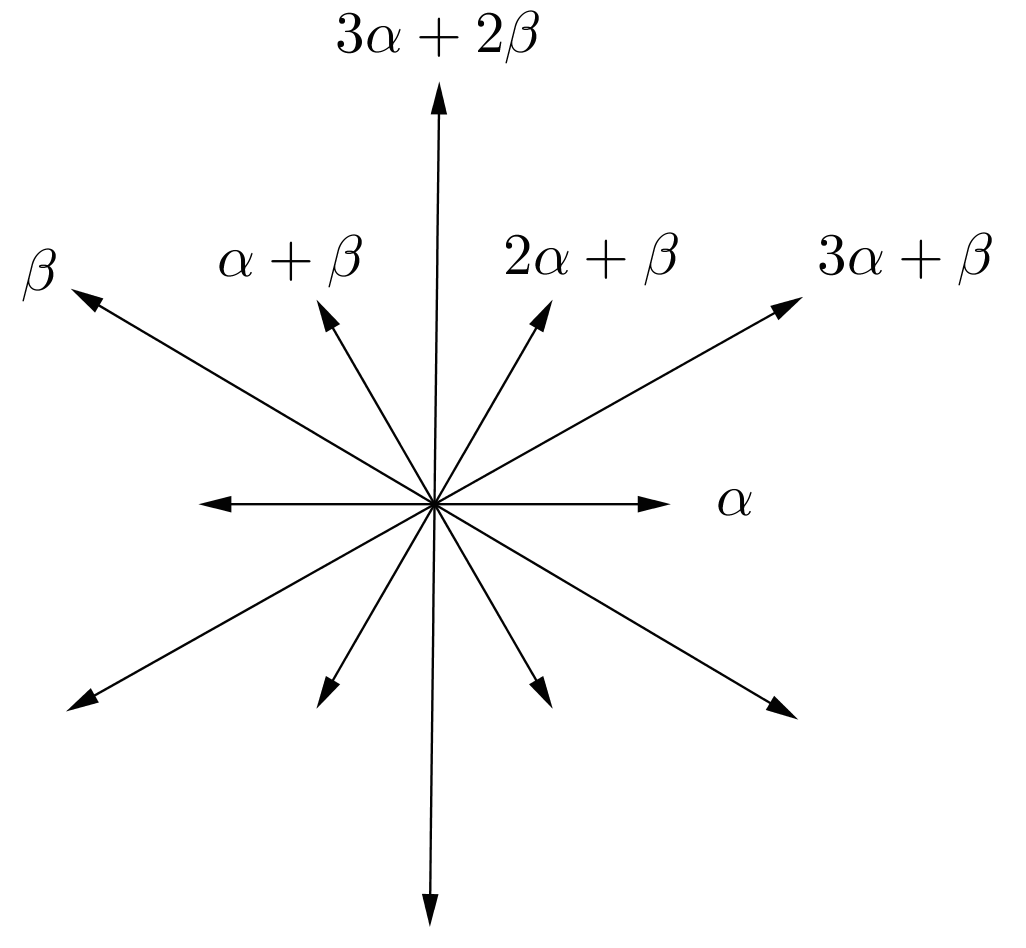}
\caption{The root system of $G_2$}\label{jhg}
\end{figure}

The exceptional Lie algebra $G_2$ is the $14$-dimensional simple Lie algebra of rank two. Its root system is 
depicted in Figure \ref{jhg}, and it has positive roots
\begin{equation}
\Phi^+=\{\alpha, ~\beta, ~\alpha+\beta, ~2\alpha+\beta, ~3\alpha+\beta, ~3\alpha+2\beta\}.
\end{equation}

$G_2$ has a Chevalley basis
\begin{equation}
\vast\{
\arraycolsep=1.7pt\def\arraystretch{1.2}
\begin{array}{llllllll}
\displaystyle H_\alpha, &  H_\beta, \\ 
\displaystyle  X_\alpha, &  X_\beta, &X_{\alpha+\beta}, &X_{2\alpha+\beta}, &X_{3\alpha+\beta},&X_{3\alpha+2\beta},\\ 
 \displaystyle  Y_\alpha, &  Y_\beta, &Y_{\alpha+\beta}, &Y_{2\alpha+\beta}, &Y_{3\alpha+\beta},&Y_{3\alpha+2\beta}
\end{array}\vast\},
\end{equation}
with nonzero commutation relations \cite{may}, where $X_{-\mu}=Y_\mu$ for $\mu \in \Phi^+$,

\begin{equation}
\begin{array}{lllllllllll}
\displaystyle [H_\mu, X_\nu] =  \frac{2(\mu, \nu)}{(\mu, \mu)} X_\nu, \quad [X_\mu, X_{-\mu}]= H_\mu, ~ \mu, \nu \in \Phi,\\
\displaystyle [X_\mu, X_\nu]= N_{\mu, \nu} X_{\mu+\nu}, \quad \mu, \nu, \mu+\nu \in \Phi,  \\
\displaystyle 1= N_{\beta, \alpha} =N_{\beta, 3\alpha+\beta}=N_{3\alpha+\beta, -(3\alpha+2\beta)} =N_{2\alpha+\beta,-(3\alpha+\beta) }=\\
N_{2\alpha+\beta, -(3\alpha+2\beta)} = N_{-(3\alpha+2\beta), \alpha+\beta}=N_{-(3\alpha+2\beta), \beta} =N_{-(3\alpha+\beta), \alpha}=\\
\displaystyle N_{-(\alpha+\beta),\beta},\\
2=N_{\alpha+\beta, \alpha}=N_{\alpha, -(2\alpha+\beta)}=N_{-(2\alpha+\beta), \alpha+\beta},\\
3=N_{\alpha, 2\alpha+\beta}=N_{\alpha,-(\alpha+\beta)}=N_{\alpha+\beta,2\alpha+\beta},\\
N_{\mu, \nu}=-N_{\nu, \mu}=-N_{-\mu, -\nu},
\end{array}
\end{equation}
and we set $(\alpha, \alpha)=1$, $(\beta, \beta)=3$, and $(\alpha, \beta)=-3/2$.

$H_\alpha$ and $H_\beta$ are the coroots; $X_\alpha,   X_\beta, X_{\alpha+\beta}, X_{2\alpha+\beta},$ $X_{3\alpha+\beta}$, and $X_{3\alpha+2\beta}$ are the positive root vectors; and 
$Y_\alpha,   Y_\beta, Y_{\alpha+\beta}, Y_{2\alpha+\beta},$ $Y_{3\alpha+\beta}$, and $Y_{3\alpha+2\beta}$ are the negative root vectors.

We highlight three distinguished subalgebras: A Cartan subalgebra $\mathfrak{h}=\langle H_\alpha, H_\beta \rangle$; a Borel subalgebra $\mathfrak{b}=\langle H_\alpha, H_\beta, X_\alpha,   X_\beta, X_{\alpha+\beta}, X_{2\alpha+\beta},$ $X_{3\alpha+\beta}, X_{3\alpha+2\beta} \rangle$; and the nilradical of the Borel subalgebra $\mathfrak{n}=\langle X_\alpha,   X_\beta,$ $X_{\alpha+\beta}, X_{2\alpha+\beta},$ $X_{3\alpha+\beta}, X_{3\alpha+2\beta} \rangle$.


The classification and identification of the subalgebras of $G_2$ are described in the following tables.  The classification is from \cite{may}, and the identification was 
established as part of the present article. The classification is organized by dimension, and further divided into regular and non-regular subalgebras. Tables \ref{g2tab0}-\ref{g2tab1} contain the classification and identification of the solvable subalgebras, Table \ref{g2tab2} contains the classification of semisimple subalgebras, and Table \ref{g2tab3} contains the classification and identification of the Levi decomposable subalgebras.

\renewcommand{\arraystretch}{1.5}
\begin{table}[h!]
\begin{tabular}{|c|c|c|c|c|c|}
  \hline 
 Representative &   Isomorphism Class\\ 
  \hline
   \multicolumn{2}{|c|}{Regular Subalgebras}  \\
   \hline
 $\langle X_\alpha \rangle$ &   $\mathfrak{n}_{1,1}$\\ 
\hline
 $\langle X_\beta \rangle$ &   $\mathfrak{n}_{1,1}$\\ 
\hline
    $\langle H_{c\alpha+d\beta}\rangle$ &   $\mathfrak{n}_{1,1}$\\ 
$(2c-3d)(-c+2d)\geq 0, (c,d) \neq (0,0)$ &  \\
\hline
 \multicolumn{2}{|c|}{Non-Regular Subalgebras}  \\
   \hline
 $\langle X_\alpha+X_{3\alpha+2\beta} \rangle$ &   $\mathfrak{n}_{1,1}$\\ 
\hline
$\langle X_\alpha+X_{\beta} \rangle$ &   $\mathfrak{n}_{1,1}$\\ 
\hline
 $\langle H_{3\alpha+2\beta}+X_\alpha \rangle$ &   $\mathfrak{n}_{1,1}$\\ 
\hline
 $\langle H_{2\alpha+\beta}+X_\beta \rangle$ &   $\mathfrak{n}_{1,1}$\\ 
\hline
 \end{tabular}
\vspace{1.5mm}
\caption{One-dimensional (solvable) subalgebras of $G_2$}\label{g2tab0}
\end{table}

\renewcommand{\arraystretch}{1.5}
\begin{table}[h!]
\begin{tabular}{|c|c|c|c|c|c|}
  \hline 
  Representative &   Isomorphism Class & Conditions\\ 
  \hline
   \multicolumn{3}{|c|}{Regular Subalgebras}  \\
   \hline
$\mathfrak{h}$ &  $2\mathfrak{n}_{1,1}$ &\\ 
\hline
$\langle H_{c\alpha+d\beta}, X_\alpha \rangle$ &  $2\mathfrak{n}_{1,1}$ & if $2c=3d$\\ 
$(2c-3d)d\geq 0$, $(c,d)\neq (0,0)$&  $\mathfrak{s}_{2,1}$ &otherwise\\
\hline
$\langle H_{c\alpha+d\beta}, X_\beta \rangle$ & $2\mathfrak{n}_{1,1}$    &if  $c=2d$\\ 
$(2d-c)c\geq 0$, $(c,d)\neq (0,0)$ & $\mathfrak{s}_{2,1}$ &otherwise \\
\hline
$\langle X_\beta,  X_{3\alpha+2\beta}\rangle$ &  $2\mathfrak{n}_{1,1}$& \\ 
\hline
$\langle X_\alpha,  X_{3\alpha+\beta}\rangle$ &  $2\mathfrak{n}_{1,1}$& \\ 
\hline
$\langle X_\alpha,  X_{3\alpha+2\beta}\rangle$ &  $2\mathfrak{n}_{1,1}$& \\ 
\hline
 \multicolumn{3}{|c|}{Non-Regular Subalgebras}  \\
   \hline
 $\langle    X_\alpha+X_\beta,  X_{3\alpha+2\beta}   \rangle$ &  $2\mathfrak{n}_{1,1}$& \\ 
\hline
 $\langle    X_\beta+X_{3\alpha+\beta},  X_{2\alpha+\beta}   \rangle$ &  $2\mathfrak{n}_{1,1}$& \\ 
\hline
 $\langle    X_{\alpha+\beta}+X_{3\alpha+\beta},  X_{3\alpha+2\beta}   \rangle$ &  $2\mathfrak{n}_{1,1}$ &\\ 
\hline
 $\langle  H_{3\alpha+\beta},   X_\alpha+X_{3\alpha+2\beta}  \rangle$ &   $\mathfrak{s}_{2,1}$ &\\ 
\hline
 $\langle  H_{9\alpha+5\beta},  X_\alpha+X_{\beta}  \rangle$ &   $\mathfrak{s}_{2,1}$ &\\ 
\hline
$\langle  H_{3\alpha+2\beta}+X_\alpha,   X_{3\alpha+\beta}  \rangle$ & $\mathfrak{s}_{2,1}$  & \\ 
\hline
 $\langle  H_{3\alpha+2\beta}+X_\alpha,   X_{3\alpha+2\beta}  \rangle$ &   $\mathfrak{s}_{2,1}$&\\ 
\hline
 $\langle  H_{2\alpha+\beta}+X_\beta,    X_{\alpha+\beta}  \rangle$ &  $\mathfrak{s}_{2,1}$& \\ 
\hline
 $\langle  H_{2\alpha+\beta}+X_\beta,   X_{2\alpha+\beta}  \rangle$ &  $\mathfrak{s}_{2,1}$&  \\ 
\hline
 $\langle  H_{2\alpha+\beta}+X_\beta,  X_{3\alpha+2\beta}  \rangle$ &   $\mathfrak{s}_{2,1}$&\\ 
\hline
 \end{tabular}
\vspace{1.5mm}
\caption{Two-dimensional (solvable) subalgebras of $G_2$}
\end{table}

\renewcommand{\arraystretch}{1.5}
\begin{table}[h!]
\begin{tabular}{|c|c|c|c|c|c|}
  \hline 
 Representative &   Isomorphism Class & Conditions\\ 
   \hline
 $\langle  H_{3\alpha+2\beta}, H_\alpha, X_\alpha\rangle$ & $\mathfrak{n}_{1, 1} \oplus \mathfrak{s}_{2, 1}$ & \\ 
\hline
 $\langle H_{2\alpha+\beta}, H_\beta, X_\beta\rangle$ & $\mathfrak{n}_{1, 1} \oplus \mathfrak{s}_{2, 1}$& \\ 
\hline
  $\langle  H_{c\alpha+d\beta},  X_\alpha, X_{3\alpha+\beta} \rangle$ &  $\mathfrak{n}_{1,1} \oplus \mathfrak{s}_{2,1}$ & if $2c=3d$ or $c=d$ \\ 
 &  $\mathfrak{s}_{3, 1, A}$  & otherwise\\
 & &$A= \frac{2c-3d}{3(c-d)}$, or $ \frac{3(c-d)}{2c-3d}$\\
& &$|A| =\min \big\{ \big| \frac{2c-3d}{3(c-d)}\big|, \big|\frac{3(c-d)}{2c-3d}\big| \big\}$ \\
&&$\arg(A) \leq \pi$ if $|A|=1$ \\
\hline
   $\langle  H_{c\alpha+d\beta},  X_\alpha, X_{3\alpha+2\beta} \rangle$  &  $\mathfrak{n}_{1,1} \oplus \mathfrak{s}_{2,1}$ & if $ 2c=3d$ or $d=0$\\ 
 & $\mathfrak{s}_{3, 1, A}$   &otherwise\\
  & &$A= \frac{2c-3d}{3d}$, or $\frac{3d}{2c-3d}$\\
& &$|A| =\min \big\{ \big|\frac{2c-3d}{3d}\big|, \big|\frac{3d}{2c-3d}\big| \big\}$ \\
&  &$\arg(A) \leq \pi$ if $|A|=1$ \\
\hline
$\langle    H_{c\alpha+d\beta}, X_\beta, X_{3\alpha+2\beta} \rangle$&   $\mathfrak{n}_{1,1} \oplus \mathfrak{s}_{2,1}$  &if  $c=2d$\\ 
$(3d-c)(c-d) \geq 0, (c,d)\neq (0,0)$& $\mathfrak{s}_{3, 1, A}$  & otherwise\\
&&$A=  \frac{c-2d}{ d}$, or $ \frac{d}{c-2d}$\\
&& $|A| =\min \big\{ \big| \frac{c-2d}{ d}\big|, \big| \frac{d}{c-2d}\big| \big\}$ \\
&& $\arg(A) \leq \pi$ if $|A|=1$\\
\hline
 $\langle X_\alpha, X_{2\alpha+\beta}, X_{3\alpha+\beta}\rangle$ &  $\mathfrak{n}_{3,1}$ &\\ 
\hline
 $\langle X_\beta, X_{\alpha+\beta}, X_{3\alpha+2\beta}\rangle$ &   $3\mathfrak{n}_{1,1}$&\\ 
\hline
 $\langle X_\beta, X_{3\alpha+\beta}, X_{3\alpha+2\beta}\rangle$ & $\mathfrak{n}_{3,1}$  &\\ 
\hline
$\langle X_\alpha, X_{3\alpha+\beta}, X_{3\alpha+2\beta}\rangle$ &  $3\mathfrak{n}_{1,1}$&\\ 
\hline
 \end{tabular}
\vspace{1.5mm}
\caption{Regular three-dimensional solvable subalgebras of $G_2$}
\end{table}

\renewcommand{\arraystretch}{1.5}
\begin{table}[h!]
\begin{tabular}{|c|c|c|c|c|c|}
  \hline 
 Representative &   Isomorphism Class \\
   \hline
 $\langle X_\alpha+X_\beta, X_{3\alpha+\beta}, X_{3\alpha+2\beta}\rangle$& $\mathfrak{n}_{3,1}$  \\ 
\hline
  $\langle X_\beta+X_{3\alpha+\beta}, X_{2\alpha+\beta}, X_{3\alpha+2\beta}\rangle$&   $3\mathfrak{n}_{1,1}$ \\ 
\hline
  $\langle X_{\alpha+\beta}+X_{3\alpha+\beta}, X_{2\alpha+\beta}, X_{3\alpha+2\beta}\rangle$& $\mathfrak{n}_{3,1}$   \\ 
\hline
  $\langle X_\beta+X_{\alpha+\beta}+\lambda X_{3\alpha+\beta}, X_{2\alpha+\beta}, X_{3\alpha+2\beta}\rangle$& $\mathfrak{n}_{3,1}$   \\ 
  $\lambda \in \mathbb{C}^*$ &\\
\hline
  $\langle H_{9\alpha+5\beta}, X_{\alpha}+X_\beta, X_{3\alpha+2\beta}\rangle$&  $\mathfrak{s}_{3, 1, A=1/5}$  \\ 
\hline
  $\langle H_{3\alpha+2\beta}, X_{\beta}+X_{3\alpha+\beta}, X_{2\alpha+\beta}\rangle$&  $\mathfrak{s}_{3, 1, A=1}$ \\ 
\hline
  $\langle H_{3\alpha+2\beta}, X_{\alpha+\beta}+X_{3\alpha+\beta}, X_{3\alpha+2\beta}\rangle$&  $\mathfrak{s}_{3, 1, A=1/2}$ \\ 
\hline
  $\langle H_{\alpha}+X_{3\alpha+2\beta}, X_{\alpha}, X_{3\alpha+\beta}\rangle$&  $\mathfrak{s}_{3, 1, A=2/3}$   \\ 
\hline
  $\langle H_{\alpha}+X_{3\alpha+2\beta}, X_{\beta}, X_{\alpha+\beta}\rangle$& $\mathfrak{s}_{3, 1, A=1/3}$   \\ 
\hline
  $\langle H_{\alpha+\beta}+X_{3\alpha+\beta}, X_{\alpha}, X_{3\alpha+2\beta}\rangle$& $\mathfrak{s}_{3, 1, A=-1/3}$   \\ 
\hline
  $\langle H_{\alpha+\beta}+X_{3\alpha+\beta}, X_{\beta}, X_{3\alpha+2\beta}\rangle$& $\mathfrak{s}_{3, 2}$   \\ 
\hline
 $\langle H_{\beta}+X_{2\alpha+\beta}, X_{\alpha}, X_{3\alpha+\beta}\rangle$&  $\mathfrak{s}_{3, 2}$ \\ 
\hline
  $\langle H_{\beta}+X_{2\alpha+\beta}, X_{\beta}, X_{3\alpha+2\beta}\rangle$&  $\mathfrak{s}_{3, 1, A=1/2}$  \\ 
\hline
 $\langle H_{3\alpha+\beta}+X_{\alpha+\beta}, X_{\beta}, X_{3\alpha+2\beta}\rangle$&  $\mathfrak{s}_{3, 1, A=-1}$  \\ 
\hline
 \end{tabular}
\vspace{1.5mm}
\caption{Non-regular three-dimensional solvable subalgebras of $G_2$}
\end{table}

\begin{footnotesize}
\renewcommand{\arraystretch}{1.5}
\begin{table}[h!]
\begin{tabular}{|c|c|c|c|c|c|}
  \hline 
 Representative &   Isomorphism Class & Conditions\\ 
   \hline
 $\langle H_\alpha, H_\beta, X_\beta, X_{3\alpha+2\beta}\rangle$ &  $2\mathfrak{s}_{2, 1}$  & \\ 
\hline
 $\langle H_\alpha, H_\beta, X_\alpha, X_{3\alpha+\beta}  \rangle$ &  $2\mathfrak{s}_{2, 1}$  & \\ 
\hline
 $\langle H_\alpha, H_\beta, X_\alpha, X_{3\alpha+2\beta}   \rangle$ &  $2\mathfrak{s}_{2, 1}$  & \\ 
\hline
  $\langle  X_\alpha, X_{2\alpha+\beta}, X_{3\alpha+\beta}, X_{3\alpha+2\beta}   \rangle$& $\mathfrak{n}_{3,1}\oplus \mathfrak{n}_{1,1}$& \\
\hline
 $\langle  X_\beta, X_{2\alpha+\beta}, X_{3\alpha+\beta}, X_{3\alpha+2\beta}   \rangle$&$\mathfrak{n}_{3,1}\oplus \mathfrak{n}_{1,1}$& \\
\hline
   $ \langle  H_{c\alpha+d\beta}, X_\alpha, X_{2\alpha+\beta}, X_{3\alpha+\beta}\rangle$ &  $\mathfrak{s}_{4, 11}$&  if $c=0$ or $2c=3d$ \\ 
& $\mathfrak{s}_{4, 6}$ & if $c=d$ \\
&$\mathfrak{s}_{4, 8, A}$ & otherwise \\
&& $A= \frac{2c-3d}{c}$ or $\frac{c}{2c-3d}$\\
&& $|A|= \min \big\{ \big|\frac{2c-3d}{c}\big|, \big| \frac{c}{2c-3d}\big| \big\}$\\
&& $\arg (A) \leq \pi$ if $|A|=1$ \\
\hline
   $\langle H_{c\alpha+d\beta}, X_\beta, X_{\alpha+\beta}, X_{3\alpha+2\beta}\rangle$ & $\mathfrak{n}_{1,1}\oplus\mathfrak{s}_{3, 1, A=1/3}$ & if $ c=2d$ or $d=0$  \\ 
& $\mathfrak{n}_{1,1}\oplus\mathfrak{s}_{3, 1, A=-1}$  &if $c=3d$ \\
&   $\mathfrak{s}_{4, 3, A , B}$& otherwise\\
&&  $A,B \in \big\{ \frac{-c+3d}{-3c+3d}, \frac{d}{-c+d},  \frac{-3c+6d}{-c+3d},$\\
&&  $ \frac{3d}{-c+3d},  \frac{-c+2d}{d},  \frac{-c+3d}{3d}\big\}$ \\
&& $0<|B| \leq |A| \leq 1$\\
\hline
   $\langle H_{c\alpha+d\beta}, X_\beta, X_{3\alpha+\beta}, X_{3\alpha+2\beta}\rangle$ &    $\mathfrak{s}_{4, 11}$ & if $c=2d$  \\ 
$(2c-3d)d \geq 0, (c,d)\neq (0,0)$&  $\mathfrak{s}_{4, 6}$&  if $d=0$\\
&$\mathfrak{s}_{4, 8, A}$& otherwise\\
&& $A=  \frac{c-d}{-c+2d}$ or $\frac{-c+2d}{c-d}$\\
&& $|A|= \min \big\{  \big|\frac{c-d}{-c+2d}\big|,  \big|\frac{-c+2d}{c-d}\big|\big\}$\\
&& $\arg (A) \leq \pi$ if $|A|=1$ \\
\hline
   $\langle H_{c\alpha+d\beta}, X_\alpha, X_{3\alpha+\beta}, X_{3\alpha+2\beta}\rangle$ & $\mathfrak{n}_{1,1}\oplus\mathfrak{s}_{3, 1, A=1/2}$&  if $2c=3d$ \\ 
  &  $\mathfrak{n}_{1,1}\oplus\mathfrak{s}_{3, 1, A=-1/3}$ & if $ c=d$  \\
&  $\mathfrak{n}_{1,1}\oplus\mathfrak{s}_{3, 1, A=2/3}$&  if $d=0$ \\
&   $\mathfrak{s}_{4, 3, A , B}$ &otherwise\\
& &$A, B \in \big\{ \frac{2c-3d}{3c-3d}, \frac{d}{c-d},  \frac{3c-3d}{2c-3d}, $\\
& & $ \frac{3d}{2c-3d},  \frac{c-d}{d},  \frac{2c-3d}{3d}\big\}$ \\
& &$0<|B| \leq |A| \leq 1$\\
\hline
 \end{tabular}
\vspace{1.5mm}
\caption{Regular four-dimensional solvable subalgebras of $G_2$}
\end{table}
\end{footnotesize}

\renewcommand{\arraystretch}{1.5}
\begin{table}[h!]
\begin{tabular}{|c|c|c|c|c|c|}
  \hline 
 Representative &   Isomorphism Class \\ 
  \hline
 $\langle X_\alpha+X_\beta, X_{2\alpha+\beta}, X_{3\alpha+\beta}, X_{3\alpha+2\beta}\rangle$ & $\mathfrak{n}_{4, 1}$ \\ 
\hline
  $\langle X_\beta+X_{2\alpha+\beta}, X_{\alpha+\beta}, X_{3\alpha+\beta}, X_{3\alpha+2\beta}\rangle$&  $\mathfrak{n}_{3, 1}\oplus \mathfrak{n}_{1,1}$ \\ 
\hline
 $\langle H_{9\alpha+5\beta}, X_{\alpha}+X_{\beta}, X_{3\alpha+\beta}, X_{3\alpha+2\beta}\rangle$ & $\mathfrak{s}_{4, 8, A=1/4}$   \\ 
\hline
 $\langle H_{3\alpha+2\beta}, X_{\beta}+X_{3\alpha+\beta}, X_{2\alpha+\beta}, X_{3\alpha+2\beta}\rangle$ & $\mathfrak{s}_{4, 3, A=B=1/2}$    \\ 
\hline
 $\langle H_{3\alpha+2\beta}, X_{\alpha+\beta}+X_{3\alpha+\beta}, X_{2\alpha+\beta}, X_{3\alpha+2\beta}\rangle$ & $\mathfrak{s}_{4, 8, A=1}$  \\ 
\hline
  $\langle H_{3\alpha+2\beta}, X_{\beta}+X_{\alpha+\beta}+\lambda X_{3\alpha+\beta}, X_{2\alpha+\beta}, X_{3\alpha+2\beta}\rangle$, $\lambda \in \mathbb{C}^*$& $\mathfrak{s}_{4, 8, A=1}$  \\ 
\hline
$\langle H_{2\alpha+\beta}+X_\beta, X_{\alpha+ \beta}, X_{3\alpha+\beta}, X_{3\alpha+2\beta}\rangle$  &  $\mathfrak{s}_{4, 4, A=1/3}$ \\ 
\hline
 $\langle H_{3\alpha+\beta}+X_{\alpha+\beta}, X_{\beta}, X_{3\alpha+\beta}, X_{3\alpha+2\beta}\rangle$ &  $\mathfrak{s}_{4, 8, A=-1/2}$  \\ 
\hline
 $\langle H_{\beta}+X_{2\alpha+\beta}, X_{\alpha}, X_{3\alpha+\beta}, X_{3\alpha+2\beta}\rangle$ & $\mathfrak{s}_{4, 4, A=-1}$  \\ 
\hline
$\langle H_{\beta}+X_{2\alpha+\beta}, X_{\beta}, X_{\alpha+\beta}, X_{3\alpha+2\beta}\rangle$  &  $\mathfrak{s}_{4, 4, A=2}$ \\ 
\hline
  $\langle H_{\alpha+\beta}+X_{3\alpha+\beta}, X_{\beta}, X_{\alpha+\beta}, X_{3\alpha+2\beta}\rangle$&  $\mathfrak{s}_{4, 4, A=2/3}$  \\ 
\hline
  $\langle H_{\alpha}+X_{3\alpha+2\beta}, X_{\alpha}, X_{2\alpha+\beta}, X_{3\alpha+\beta}\rangle$&   $\mathfrak{s}_{4, 8, A=1/2}$ \\ 
\hline
 \end{tabular}
\vspace{1.5mm}
\caption{Non-regular four-dimensional solvable subalgebras of $G_2$}
\end{table}

\begin{scriptsize}
\renewcommand{\arraystretch}{1.5}
\begin{table}[h!]
\begin{tabular}{|c|c|c|c|c|c|}
  \hline 
 Representative &   Isomorphism Class & Conditions\\ 
  \hline
   \multicolumn{3}{|c|}{Regular Subalgebras}  \\
   \hline
 $\langle H_\alpha, H_\beta, X_\alpha, X_{2\alpha+\beta}, X_{3\alpha+\beta}\rangle$ & $\mathfrak{s}_{5, 44}$&  \\ 
\hline
 $\langle H_\alpha, H_\beta, X_\beta, X_{\alpha+\beta}, X_{3\alpha+2\beta}\rangle$ & $\mathfrak{s}_{5, 41, A=B=1/3}$ & \\ 
\hline
 $\langle H_\alpha, H_\beta, X_\beta, X_{3\alpha+\beta}, X_{3\alpha+2\beta}\rangle$ &   $\mathfrak{s}_{5, 44}$&\\ 
\hline
  $\langle H_\alpha, H_\beta, X_\alpha, X_{3\alpha+\beta}, X_{3\alpha+2\beta}\rangle$ &   $\mathfrak{s}_{5, 41, A=2/3, B=-1/3}$ &\\ 
\hline
  $\langle  H_{c\alpha+d\beta}, X_\alpha, X_{2\alpha+\beta}, X_{3\alpha+\beta}, X_{3\alpha+2\beta}\rangle$ & $\mathfrak{s}_{5,30, A=2}$& if $2c=3d$ \\ 
 & $\mathfrak{s}_{5,30, A=-1}$& if $ c=0$ \\
& $\mathfrak{s}_{5,17, A=3}$& if $ c=d$\\
&  $\mathfrak{n}_{1,1}\oplus \mathfrak{s}_{4,8, A=1/2}$& if $d=0$\\
& $\mathfrak{s}_{5, 22, A, B}$& otherwise\\
&&  $A= \frac{2c-3d}{c}$ or $\frac{c}{2c-3d}$\\
&& $|A|=\min\big\{  \big| \frac{2c-3d}{c}\big|, \big| \frac{c}{2c-3d}\big| \big\}$\\
&& $\arg (A) \leq \pi$ if $|A|=1$ \\
&&  $B=  \frac{3d}{c}$ or $\frac{3d}{2c-3d}$ \\
\hline
  $\langle   H_{c\alpha+d\beta}, X_\beta, X_{2\alpha+\beta}, X_{3\alpha+\beta}, X_{3\alpha+2\beta}\rangle$ &  $\mathfrak{s}_{5,30, A=2/3}$& if $c=2d$ \\ 
 & $\mathfrak{s}_{5,30, A=1/3}$& if $c=d$\\
& $\mathfrak{s}_{5,17, A=1/3}$& if $d=0$\\
&  $\mathfrak{n}_{1,1}\oplus \mathfrak{s}_{4,8, A=-1/2}$& if $c=0$\\
& $\mathfrak{s}_{5, 22, A, B}$& otherwise\\
&&  $A=\frac{c-d}{-c+2d}$ or $ \frac{-c+2d}{c-d}$\\
&& $|A|=\min\big\{  \big|\frac{c-d}{-c+2d}\big|, \big| \frac{-c+2d}{c-d}\big| \big\}$\\
&& $\arg (A) \leq \pi$ if $|A|=1$ \\
&&  $B= \frac{c}{3c-3d}$ or $\frac{c}{-3c+6d}$ \\
\hline
 $\langle X_\alpha, X_{\alpha+\beta}, X_{2\alpha+\beta}, X_{3\alpha+\beta}, X_{3\alpha+2\beta}\rangle$ &  $\mathfrak{n}_{5, 2}$ & \\ 
\hline
 $\langle X_\beta, X_{\alpha+\beta}, X_{2\alpha+\beta}, X_{3\alpha+\beta}, X_{3\alpha+2\beta}\rangle$ &  $\mathfrak{n}_{5, 3}$ & \\ 
\hline
 \multicolumn{3}{|c|}{Non-Regular Subalgebras}  \\
   \hline
 $\langle X_\alpha+X_\beta, X_{\alpha+\beta}, X_{2\alpha+\beta}, X_{3\alpha+\beta}, X_{3\alpha+2\beta}\rangle$& $\mathfrak{n}_{5,6}$& \\ 
\hline
 $\langle H_{9\alpha+5\beta}, X_\alpha+X_{\beta}, X_{2\alpha+\beta},
X_{3\alpha+\beta}, X_{3\alpha+2\beta}\rangle$& $\mathfrak{s}_{5, 35, A=3}$ & \\ 
\hline
 $\langle H_{3\alpha+2\beta}, X_\beta+X_{2\alpha+\beta}, X_{\alpha+\beta},
X_{3\alpha+\beta}, X_{3\alpha+2\beta}\rangle$& $\mathfrak{s}_{5, 22, A=B=1}$& \\ 
\hline
 $\langle H_{3\alpha+\beta}+X_{\alpha+\beta}, X_\beta, X_{2\alpha+\beta},
X_{3\alpha+\beta}, X_{3\alpha+2\beta}\rangle$&$\mathfrak{s}_{5, 26, A=-1/2}$ & \\ 
\hline
 $\langle H_{3\alpha+\beta}+X_{\alpha+\beta}, X_\alpha, X_{2\alpha+\beta},
X_{3\alpha+\beta}, X_{3\alpha+2\beta}\rangle$&  $\mathfrak{s}_{5, 21}$&\\ 
\hline
 $\langle H_{\alpha+\beta}+X_{3\alpha+\beta}, X_\beta, X_{\alpha+\beta},
X_{2\alpha+\beta}, X_{3\alpha+2\beta}\rangle$& $\mathfrak{s}_{5, 26, A=1/2}$& \\ 
\hline
 \end{tabular}
\vspace{1.5mm}
\caption{Five-dimensional solvable subalgebras of $G_2$}
\end{table}
\end{scriptsize}

\begin{landscape}
\begin{scriptsize}
\renewcommand{\arraystretch}{1.5}
\begin{table}[h!]
\begin{tabular}{|c|c|c|c|c|c|}
  \hline 
 Representative &   Isomorphism Class & Conditions\\ 
  \hline
   \multicolumn{3}{|c|}{Regular Subalgebras}  \\
   \hline
  $\langle H_{c\alpha+d\beta }, X_\alpha, X_{\alpha+\beta}, X_{2\alpha+\beta}, X_{3\alpha+\beta}, X_{3\alpha+2\beta}\rangle$ & $\mathfrak{s}_{6, 150}$& if  $c=d$\\ 
$(2d-c)c \geq 0, (c,d)\neq (0,0)$& $\mathfrak{s}_{6, 156}$& if $2c=3d$\\
& $\mathfrak{s}_{6, 155, A=-1}$& if $c=0$\\
& $\mathfrak{s}_{6, 155, A}$& otherwise\\
&& $A= \frac{2c-3d}{-c+3d}$ or $\frac{-c+3d}{2c-3d}$\\
&& $|A| =\min \big\{ \big|\frac{2c-3d}{-c+3d}\big|, \big|\frac{-c+3d}{2c-3d}\big| \big\}$\\
&& $\arg (A) \leq \pi$ if $|A|=1$ \\
\hline
  $\langle H_{c\alpha+d\beta}, X_\beta, X_{\alpha+\beta}, X_{2\alpha+\beta}, X_{3\alpha+\beta}, X_{3\alpha+2\beta} \rangle$ & $\mathfrak{s}_{6, 162, A=1/3}$& if  $d=0$ \\ 
$(2c-3d)d\geq 0, (c,d)\neq (0,0)$ & $\mathfrak{s}_{6, 178, A=1/3}$& if   $c=2d$\\
& $\mathfrak{s}_{6, 178, A=3}$& if  $c=3d$\\
&    $\mathfrak{s}_{6, 168, A, B}$& otherwise\\
&& $A, B \in \big\{ \pm \frac{-2c+3d}{d}, \pm \frac{-2c+3d}{3d} \big\}$\\
&&  $|B|\leq |A|$, $\arg(A), \arg(B) <\pi$\\
&& if $|A|=|B|$ then $\arg(B)\leq \arg(A)$ \\ 
\hline
 $\langle  H_\alpha, H_\beta, X_\alpha, X_{2\alpha+\beta}, X_{3\alpha+\beta}, X_{3\alpha+2\beta}\rangle$ &$\mathfrak{s}_{6,234,A=-1, B=2}$& \\
\hline 
 $\langle  H_\alpha, H_\beta, X_\beta, X_{2\alpha+\beta}, X_{3\alpha+\beta}, X_{3\alpha+2\beta}\rangle$  &$\mathfrak{s}_{6,234,A=1/3, B=2/3}$& \\
\hline
 $\mathfrak{n}$ &  $\mathfrak{n}_{6,18}$& \\ 
\hline
 \multicolumn{3}{|c|}{Non-Regular Subalgebras}  \\
   \hline
 $\langle H_{9\alpha+5\beta}, X_\alpha+ X_\beta, X_{\alpha+\beta}, X_{2\alpha+\beta}, X_{3\alpha+\beta}, X_{3\alpha+2\beta}\rangle$& $\mathfrak{s}_{6, 197}$ &\\ 
\hline
 $\langle H_{3\alpha+2\beta}+X_\alpha, X_\beta, X_{\alpha+\beta}, X_{2\alpha+\beta}, X_{3\alpha+\beta}, X_{3\alpha+2\beta}\rangle$ & $\mathfrak{s}_{6, 171}$  &\\ 
\hline
 $\langle H_{2\alpha+\beta}+X_\beta, X_\alpha, X_{\alpha+\beta}, X_{2\alpha+\beta}, X_{3\alpha+\beta}, X_{3\alpha+2\beta}\rangle$ &$\mathfrak{s}_{6, 153}$ &  \\ 
\hline
 \end{tabular}
\vspace{1.5mm}
\caption{Six-dimensional solvable subalgebras of $G_2$}
\end{table}
\end{scriptsize}
\end{landscape}

\begin{footnotesize}
\renewcommand{\arraystretch}{1.5}
\begin{table}[h!]
\begin{tabular}{|c|c|c|c|c|c|}
  \hline 
 Dimension & Representative &   Isomorphism Class & Conditions\\ 
  \hline
   \multicolumn{4}{|c|}{Regular Subalgebras}  \\
   \hline
$7$  & $\mathfrak{h}+\langle  X_\alpha, X_{\alpha+\beta}, X_{2\alpha+\beta}, X_{3\alpha+\beta}, X_{3\alpha+2\beta}\rangle$ & $\mathfrak{s}_{7, 1}$&  \\ 
\hline
$7$  & $\mathfrak{h}+\langle  X_\beta, X_{\alpha+\beta}, X_{2\alpha+\beta}, X_{3\alpha+\beta}, X_{3\alpha+2\beta}\rangle$ &  $\mathfrak{s}_{7, 2}$& \\ 
\hline
$7$  &  $\mathfrak{n}+\langle H_{c\alpha+d\beta} \rangle$ & $\mathfrak{s}_{7, 3,  c, 3d}$& $\mathfrak{s}_{7, 3, c, 3d} \cong \mathfrak{s}_{7, 3, c', 3d'}$\\
&& & iff $cd'=c'd$  \\ 
\hline
$8$  & $\mathfrak{b}$ &  $\mathfrak{s}_{8, 1}$ &\\ 
\hline
 \end{tabular}
\vspace{1.5mm}
\caption{Seven- and eight-dimensional solvable subalgebras of $G_2$}\label{g2tab1}
\end{table}
\end{footnotesize}

\begin{small}
\renewcommand{\arraystretch}{1.5}
\begin{table}[h!]
\begin{tabular}{|c|c|c|c|c|c|}
  \hline 
 Dimension & Representative &   Isomorphism Class\\ 
  \hline
   \multicolumn{3}{|c|}{Regular Subalgebras}  \\
   \hline
3  & $\langle H_\alpha, X_\alpha, Y_{\alpha}\rangle$ & $A_1$  \\ 
\hline
3  & $\langle H_\beta, X_\beta, Y_{\beta}\rangle$ & $A_1$  \\ 
\hline
6  & $\langle H_\alpha, H_\beta, X_\alpha, X_{3\alpha+2\beta}, Y_{\alpha}, Y_{3\alpha+2\beta}\rangle$ & $A_1\times A_1$  \\ 
\hline
8 & $\langle H_\alpha, H_\beta, X_\beta, X_{3\alpha+\beta}, X_{3\alpha+2\beta}, Y_{\beta}, Y_{3\alpha+\beta}, Y_{3\alpha+2\beta}\rangle$ & $A_2$  \\ 
\hline
 \multicolumn{3}{|c|}{Non-Regular Subalgebras}  \\
   \hline
 3 &  $\langle 2H_{3\alpha+\beta}, \sqrt{2}(Y_{\beta}+X_{3\alpha+2\beta}), \sqrt{2}(X_{\beta}+Y_{3\alpha+2\beta})\rangle$&   $A_1$\\ 
\hline
 3 &  $\langle 14H_{9\alpha+5\beta}, \sqrt{6}X_{\alpha}+\sqrt{10}X_{\beta}, \sqrt{6}Y_{\alpha}+\sqrt{10}Y_{\beta} \rangle$&   $A_1$\\ 
\hline
 \end{tabular}
\vspace{1.5mm}
\caption{Semisimple subalgebras of $G_2$}\label{g2tab2}
\end{table}
\end{small}

\begin{landscape}
\begin{small}
\renewcommand{\arraystretch}{1.5}
\begin{table}[h!]
\begin{tabular}{|c|c|c|c|c|c|}
  \hline 
 Dimension & Representative &   Isomorphism Class\\ 
  \hline
   \multicolumn{3}{|c|}{Regular Subalgebras}  \\
   \hline
 $4$ & $\langle H_{\alpha}, X_{\alpha}, Y_{\alpha}\rangle \oplus \langle H_{3\alpha+2\beta}\rangle$ &  $A_1 \oplus \mathfrak{n}_{1, 1}$ \\ 
\hline
 $4$ & $\langle  H_{\beta}, X_{\beta}, Y_{\beta} \rangle \oplus \langle H_{2\alpha+\beta}\rangle$ &  $A_1 \oplus \mathfrak{n}_{1, 1}$ \\ 
\hline
 $4$ & $\langle  H_{\alpha}, X_{\alpha}, Y_{\alpha} \rangle \oplus \langle X_{3\alpha+2\beta}\rangle$ &  $A_1 \oplus \mathfrak{n}_{1, 1}$ \\ 
\hline
 $4$ & $\langle H_{\beta}, X_{\beta}, Y_{\beta}  \rangle \oplus \langle X_{2\alpha+\beta}\rangle$ &  $A_1 \oplus \mathfrak{n}_{1, 1}$ \\ 
\hline
 $5$ & $\langle H_{\beta}, X_{\beta}, Y_{\beta} \rangle \inplus \langle X_{3\alpha +\beta}, X_{3\alpha+2\beta}\rangle$ &$A_1 \inplus 2\mathfrak{n}_{1,1}$ \\ 
\hline
 $5$ & $\langle H_{\alpha}, X_{\alpha}, Y_{\alpha} \rangle \oplus \langle H_{3\alpha+2\beta}, X_{3\alpha+2\beta}\rangle$  &   $A_1 \oplus \mathfrak{s}_{2, 1}$ \\ 
\hline
 $5$ &$\langle H_{\beta}, X_{\beta}, Y_{\beta} \rangle \oplus \langle H_{2\alpha+\beta}, X_{2\alpha+\beta}\rangle$   & $A_1 \oplus \mathfrak{s}_{2, 1}$  \\ 
\hline
 $6$ &   $\langle H_{\beta}, X_{\beta}, Y_{\beta} \rangle  \inplus \langle H_{2\alpha+\beta}, X_{3\alpha+\beta}, X_{3\alpha+2\beta}\rangle$  &  $A_1 \inplus \mathfrak{s}_{3, 1, A=1}$ \\ 
\hline
 $6$ &  $(\langle H_{\beta}, X_{\beta}, Y_{\beta} \rangle \inplus  \langle X_{3\alpha+\beta}, X_{3\alpha+2\beta}\rangle) \oplus \langle X_{2\alpha+\beta}\rangle$ &  $  (A_1 \inplus 2\mathfrak{n}_{1, 1})\oplus \mathfrak{n}_{1,1}$ \\ 
\hline
 $7$ &  $\langle H_{\beta}, X_{\beta}, Y_{\beta}\rangle  \inplus  \langle H_{2\alpha+\beta}, X_{2\alpha+\beta}, X_{3\alpha+\beta}, X_{3\alpha+2\beta}\rangle$&  
  $A_1\inplus \mathfrak{s}_{4, 3, A=1, B=2/3}$  \\
\hline
 $8$ & $\langle H_{\alpha}, X_{\alpha}, Y_{\alpha} \rangle \inplus \langle X_{\beta}, X_{\alpha+\beta}, X_{2\alpha+\beta}, X_{3\alpha+\beta}, X_{3\alpha+2\beta} \rangle$ &  $A_1 \inplus \mathfrak{n}_{5, 3}$ \\ 
\hline
 $8$ &  $\langle H_{\beta}, X_{\beta}, Y_{\beta}\rangle  \inplus  \langle X_\alpha, X_{\alpha+\beta}, X_{2\alpha+\beta}, X_{3\alpha+\beta}, X_{3\alpha+2\beta}\rangle$&    $A_1 \inplus \mathfrak{n}_{5, 2}$\\ 
\hline
$9$ &  $\langle H_{\alpha}, X_{\alpha}, Y_{\alpha} \rangle \inplus \langle H_{3\alpha+2\beta}, X_\beta, X_{\alpha+\beta}, X_{3\alpha+2\beta}, 
X_{2\alpha+\beta}, X_{3\alpha+\beta} \rangle$&  $A_1 \inplus \mathfrak{s}_{6, 168, A=B=0}$ \\
\hline
$9$ &  $\langle H_{\beta}, X_{\beta}, Y_{\beta} \rangle \inplus \langle  H_{2\alpha+\beta}, X_{\alpha+\beta}, X_{3\alpha+2\beta}, X_{2\alpha+\beta}, 
X_{3\alpha+\beta}, X_{\alpha} \rangle$&  $A_1 \inplus \mathfrak{s}_{6, 155, A=1}$ \\
\hline
 \end{tabular}
\vspace{1.5mm}
\caption{Levi decomposable subalgebras of $G_2$}\label{g2tab3}
\end{table}
\end{small}
\end{landscape}

\clearpage

\section*{Acknowledgements}

The work of A.D. was partially supported by a research grant from the Professional Staff
Congress/City University of New York (Grant No. TRADA-47-36). The
work of J.R. was partially supported by the Natural Sciences and Engineering Research Council
(Grant No. 3166-09).  We also greatly thank Matthew Wu who checked our computations in the $G_2$ case.

\appendix
\addappheadtotoc


\section{Classification of  Lie algebras}

By Levi's Theorem [\cite{levi}, Chapter II, Section 2],  a complex  Lie algebra is either semisimple, solvable, or a nontrivial semidirect sum of the first two.  A Lie algebra that is a  nontrivial semidirect sum of a semisimple algebra with a solvable algebra is called a Levi decomposable algebra. Hence, the classification of Lie algebras may be divided into three cases: Semisimple Lie algebras, solvable Lie algebras, and Levi decomposable  algebras.

\subsection{Classification of semisimple Lie algebras}

Semsimple Lie algebras over $\mathbb{C}$ are direct sums of simple Lie algebras. Simple Lie algebras over $\mathbb{C}$ were classified by Wilhelm Killing \cite{killing} and \'Elie Cartan \cite{cartan}. Killing first described the classification, and Cartan provided a rigorous proof.

Under the classification, there are four families of classical Lie algebras, and five exceptional Lie algebras. The four families of classical Lie algebras are denoted $A_n$ ($n\geq 1$), $B_n$ ($n\geq 3$), $C_n$ ($n\geq 2$), and $D_n$ ($n\geq 4$). The exceptional Lie algebras are denoted $E_6, E_7$, $E_8$, $F_4$, and $G_2$. The subscript indicates the rank of the Lie algebra. A complete description of the classification and its proof can found in any standard reference on Lie algebras such as Humphreys \cite{humphreys}.

\subsection{Classification of solvable Lie algebras}\label{solvable}

A full classification of solvable Lie algebras is not known and thought to be an impossible task.  However, partial classifications of solvable Lie algebras do exist.  Two such partial classifications are that described by \v{S}nobl and Winternitz in \cite{levi} and that of de Graaf \cite{degraafa}. The classification of solvable subalgebras in this article are described, where possible, with respect to the classification described by \v{S}nobl and Winternitz.

The classification of solvable Lie algebras  presented  by  \v{S}nobl and Winternitz in \cite{levi} is up to and including dimension  $6$ and includes only   Lie algebras indecomposable over $\mathbb{R}$.  This classification is an amalgam of results from various sources (e.g., \cite{bianchi, kruch, lie, moro, patera, shaban, turko, turko2}).

We present the classification from \cite{levi} in its entirety up to and including dimension  $3$. Below we  give a 
partial description of the classification in dimensions $4$, $5$ and $6$, including just those algebras which appear in this article.
In this classification, an algebra designated with $\mathfrak{n}$ is nilpotent, and one with $\mathfrak{s}$
is solvable, but not nilpotent. The first subscript indicates the dimension, and the second index is for enumeration. So, $\mathfrak{s}_{6, 242}$  is the $242^{\text{nd}}$ six-dimensional, solvable, non-nilpotent Lie algebra in the classification.

\begin{footnotesize}
\begin{equation}\label{onesnobl}
\begin{array}{llll}
\mathfrak{n}_{1,1} & \text{The abelian Lie algebra of dimension $1$}\\
\end{array}
\end{equation}

\begin{equation}
\arraycolsep=2pt\def\arraystretch{1.5}
\begin{array}{llll}
 \mathfrak{s}_{2,1} &      [e_2, e_1]=e_1 
\end{array}
\end{equation}

\begin{equation}\label{threesnobl}
\arraycolsep=2pt\def\arraystretch{1.5}
\begin{array}{llll}
 \mathfrak{n}_{3,1} &   [e_2, e_3]=e_1 \\
 \mathfrak{s}_{3,1} &      [e_3, e_1]=e_1, [e_3, e_2]=A e_2   \\
 & 0<|A|\leq 1, ~\text{if}~|A|=1 ~\text{then}~\arg(A)\leq \pi \\
  \mathfrak{s}_{3,2} &   [e_3, e_1]=e_1, [e_3, e_2]= e_1+e_2
\end{array}
\end{equation}

\begin{equation}\label{foursnobl}
\arraycolsep=2pt\def\arraystretch{1.5}
\begin{array}{llll}
 \mathfrak{n}_{4,1} &  [e_2, e_4]=e_1, [e_3, e_4]=e_2\\
 \mathfrak{s}_{4,2} &  [e_4, e_1]=e_1, [e_4, e_2]=e_1+e_2, [e_4, e_3]=e_2+e_3 \\
   \mathfrak{s}_{4,3} & [e_4, e_1]= e_1, [e_4, e_2]=Ae_2, [e_4, e_3]=Be_3,  \\
   & 0 <|B| \leq |A| \leq 1, (A,B) \neq (-1, -1) \\
   \mathfrak{s}_{4,4} & [e_4, e_1]=e_1, [e_4, e_2]= e_1+e_2, [e_4, e_3]=Ae_3,\\
   & A\neq 0 \\
  \mathfrak{s}_{4,6} &  [e_2, e_3]=e_1, [e_4, e_2]=e_2, [e_4, e_3]=-e_3 \\
    \mathfrak{s}_{4,8} &  [e_2, e_3]=e_1, [e_4, e_1]=(1+A)e_1, [e_4, e_2]=e_2, [e_4,e_3]=Ae_3,\\
    & 0<|A|\leq 1, ~\text{if}~|A|=1 ~\text{then}~\arg(A)< \pi\\
      \mathfrak{s}_{4,10} &  [e_2, e_3]=e_1, [e_4, e_1]=2e_1, [e_4, e_2]=e_2, [e_4, e_3]=e_2+e_3 \\
        \mathfrak{s}_{4,11} &  [e_2, e_3]=e_1, [e_4, e_1]=e_1, [e_4, e_2]=e_2 \\
          \mathfrak{s}_{4,12} &  [e_3, e_1]=e_1, [e_3, e_2] = e_2, [e_4, e_1]=-e_2, [e_4, e_2]=e_1
\end{array}
\end{equation}

\begin{equation}
\arraycolsep=2pt\def\arraystretch{1.5}
\begin{array}{llll}
 \mathfrak{n}_{5,2} & [e_3, e_4]=e_2, [e_3, e_5]=e_1, [e_4, e_5]=e_3\\
  \mathfrak{n}_{5,3} & [e_2, e_4]=e_1, [e_3, e_5]=e_1\\
   \mathfrak{n}_{5,6} & [e_2, e_5]=e_1, [e_3, e_4]=e_1, [e_3, e_5]=e_2, [e_4, e_5]=e_3\\
 \mathfrak{s}_{5,17} & [e_5, e_2]=e_2, [e_5, e_3]=-e_3, [e_5, e_4]=Ae_4,\\
 &[e_2, e_3]=e_1, \\
 &  0\leq \mathrm{Re} (A), ~\text{if}~\mathrm{Re}(A)=0 ~\text{then}~ 0< \mathrm{Im}(A)\\
   \mathfrak{s}_{5,21} & [e_5,e_1]=2e_1, [e_5, e_2]=e_2+e_3, [e_5, e_3]= e_3+e_4\\
   & [e_5, e_4]=e_4,  [e_2,e_3]=e_1\\
  \mathfrak{s}_{5,22} & [e_5, e_1]=(A+1)e_1, [e_5, e_2]=e_2, [e_5, e_3]=Ae_3, [e_5, e_4]=Be_4,\\
 & [e_2,e_3]=e_1,\\
  & 0<|A| \leq 1, B\neq 0, ~\text{if}~|A|=1~\text{then}~\mathrm{arg}(A)<\pi\\
    \mathfrak{s}_{5,26} & [e_5, e_1]=(A+1)e_1, [e_5, e_2]=e_2, [e_5, e_3]=Ae_3,\\
&  [e_5, e_4]=e_1+(A+1)e_4, [e_2,e_3]=e_1,\\
&  0<|A|\leq 1, ~\text{if}~|A|=1~\text{then}~\mathrm{arg}(A)<\pi \\
   \mathfrak{s}_{5,30} & [e_5, e_1]=e_1,   [e_5, e_2]=e_2, [e_5, e_4]= Ae_4,\\
   & [e_2, e_3]=e_1,\\
   & A\neq 0\\
 \mathfrak{s}_{5,33} &  [e_2, e_4]=e_1, [e_3, e_4]=e_2, [e_5, e_2]=-e_2,\\
 &  [e_5, e_3]=-2e_3, [e_5, e_4]=e_4\\
 \mathfrak{s}_{5,35} &  [e_2, e_4]=e_1, [e_3, e_4]=e_2, [e_5, e_1]=(A+2)e_1,\\
 &  [e_5, e_2]=(A +1)e_2, [e_5, e_3]=A e_3, [e_5, e_4]= e_4,\\
 & A \neq 0, -2\\
  \mathfrak{s}_{5,36} &   [e_2, e_4]=e_1, [e_3, e_4]=e_2, [e_5, e_1]=2e_1,\\
 &  [e_5, e_2]=e_2, [e_5, e_4]=e_4\\
 \mathfrak{s}_{5,37} &  [e_2, e_4]=e_1, [e_3, e_4]=e_2, [e_5, e_1]=e_1,\\
 &  [e_5, e_2]=e_2, [e_5, e_3]=e_3\\
  \mathfrak{s}_{5,41} &  [e_4, e_1]=e_1, [e_4, e_3]=A e_3,\\
  &  [e_5, e_2]=e_2, [e_5, e_3]=B e_3,\\& 0< |B | \leq |A | \leq 1\\
   \mathfrak{s}_{5,44} & [e_2,e_3]=e_1, [e_4,e_1]=e_1, [e_4,e_2]=e_2,\\
   & [e_5,e_2]=e_2, [e_5,e_3]=-e_3
\end{array}
\end{equation}

\begin{equation}
\arraycolsep=2pt\def\arraystretch{1.5}
\begin{array}{llll}
 \mathfrak{n}_{6, 18} & [e_2, e_6]=e_1, [e_3, e_4]=e_1, [e_3, e_5]=e_2, [e_4, e_5]=e_3, [e_5, e_6]=e_4 \\
 \mathfrak{s}_{6,150} &  [e_3, e_4]=e_2, [e_3, e_5]=e_1, [e_4, e_5]=e_3, [e_6, e_2]=3e_2,\\
 & [e_6, e_3]=e_3, [e_6, e_4]=2e_4, [e_6, e_5]=-e_5\\
  \mathfrak{s}_{6, 153} &   [e_3, e_4]=e_2, [e_3, e_5]=e_1, [e_4, e_5]=e_3, [e_6, e_1]= 3e_1,\\
  & [e_6, e_2]=e_1+3e_2, [e_6, e_3]=2e_3, [e_6, e_4]=e_4+e_5, [e_6, e_5]=e_5\\
  \mathfrak{s}_{6, 155} &   [e_3, e_4]=e_2, [e_3, e_5]=e_1, [e_4, e_5]=e_3, [e_6, e_1]= (2A+1)e_1,\\
  & [e_6, e_2]= (A+2) e_2, [e_6, e_3]=(A+1)e_3, [e_6, e_4]=e_4,\\
  & [e_6, e_5]=Ae_5, \\
  & 0<|A|\leq 1, A\neq -1/2, \mathrm{arg}(A)\leq \pi~\text{if}~|A|=1\\
   \mathfrak{s}_{6, 156} &  [e_3, e_4]=e_2, [e_3, e_5]=e_1, [e_4, e_5]=e_3, [e_6, e_1]=e_1,\\
   & [e_6, e_2]= 2e_2, [e_6, e_3]=e_3, [e_6, e_4]=e_4\\
    \mathfrak{s}_{6, 162} & [e_2, e_4]=e_1, [e_3, e_5]=e_1, [e_6, e_2]=e_2, [e_6, e_3]=Ae_3\\
    & [e_6, e_4]=-e_4, [e_6, e_5]=-Ae_5,\\
    &  0<|A|\leq 1~\text{and}~ \mathrm{arg}(A)<\pi, ~\text{if}~ |A|=1 ~\text{then}~ \mathrm{arg}(A)<\pi/2\\
     \mathfrak{s}_{6, 168} &  [e_2, e_4]=e_1, [e_3, e_5]=e_1, [e_6, e_1]=2e_1, [e_6, e_2]=(A+1)e_2,\\
     & [e_6, e_3]=(B+1)e_3, [e_6, e_4]= (1-A)e_4, [e_6, e_5]= (1-B)e_5,\\
     & |B| \leq |A|,~ \mathrm{arg}(A), \mathrm{arg}(B) < \pi, A, B \neq 1. \\
     & \text{If}~ |A|=|B|~\text{then}~ \mathrm{arg}(B) \leq \mathrm{arg}(A).\\
      \mathfrak{s}_{6, 171} & [e_2, e_4]=e_1, [e_3, e_5]=e_1,  [e_6, e_1]=2e_1, [e_6, e_2]=e_2+e_3,\\
& [e_6, e_3]=e_3-e_5, [e_6, e_4]=e_4, [e_6, e_5]=-e_4+e_5\\
     \mathfrak{s}_{6, 178} & [e_2, e_4]=e_1, [e_3, e_5]=e_1,  [e_6, e_1]=2e_1, [e_6, e_2]=(1+A) e_2,\\
     & [e_6, e_3]=2e_3, [e_6, e_4]=(1-A)e_4, \\
     & A\neq 1, \mathrm{arg}(A) <\pi\\
      \mathfrak{s}_{6, 197} & [e_2, e_5]=e_1, [e_3, e_4]=e_1, [e_3, e_5]=e_2, [e_4, e_5]=e_3,\\
      &  [e_6, e_1]=5e_1, [e_6, e_2]=4e_2, [e_6, e_3]=3e_3, [e_6, e_4]=2e_4,\\
      & [e_6, e_5]=e_5\\
       \mathfrak{s}_{6, 234} & [e_2, e_3]=e_1, [e_5, e_1]=e_1, [e_5, e_2]=e_2, [e_5, e_4]=Ae_4\\
       & [e_6, e_1]=e_1, [e_6, e_3]=e_3, [e_6, e_4]=Be_4,\\
       & |A|\leq |B|, B\neq 0. ~\text{If}~|A|=|B|~\text{then}~\mathrm{arg}(A)\leq \mathrm{arg}(B)\\
 \mathfrak{s}_{6,242} &  [e_2, e_4]=e_1, [e_3, e_4]=e_2, [e_5, e_1]=2e_1, \\
  & [e_5, e_2]=e_2, [e_5, e_4]=e_4, [e_6, e_1]=e_1,\\
  & [e_6, e_2]=e_2, [e_6, e_3]=e_3
\end{array}
\end{equation}
\end{footnotesize}

We provide multiplication tables for the remaining solvable Lie algebras which appear in this article. These solvable Lie algebras 
do not appear in \cite{levi}, but we have chosen  notation consistent with that in \cite{levi}.  For these solvable Lie algebras, we also identify their nilradical.

\begin{itemize}
\item $\mathfrak{s}_{7, 1}$, nilradical  $\mathfrak{n}_{5, 2}$

\begin{center}
\begin{tabular}{ |c| c | c | c | c |c|c|c|c|c|c|c|c|}
\hline
 & $e_1$& $e_2$&$e_3$ &$e_4$ &$e_5$ & $e_6$& $e_7$\\
\hline
$e_1$ & $0$ &$0$&$0$&$0$&$0$&$-3e_1$& $e_1$\\ 
\hline
$e_2$ & $0$ &$0$&$0$&$0$&$0$&$0$&$-e_2$ \\ 
\hline
$e_3$ & $0$ &$0$&$0$&$e_2$&$e_1$&$-e_3$&$0$ \\ 
\hline
$e_4$ & $0$ &$0$&$-e_2$&$0$&$e_3$&$e_4$&$-e_4$ \\ 
\hline
$e_5$ &  $0$&$0$&$-e_1$&$-e_3$& $0$&$-2e_5$&$e_5$\\ 
\hline
$e_6$ & $3e_1$ &$0$&$e_3$&$-e_4$&$2e_5$&$0$&$0$ \\ 
\hline
$e_7$ &$-e_1$  &$e_2$&$0$&$e_4$&$-e_5$&$0$&$0$ \\ 
\hline
\end{tabular}
\end{center}

\vspace{2mm}

\item $\mathfrak{s}_{7, 2}$, nilradical  $\mathfrak{n}_{5, 3}$

\begin{center}
\begin{tabular}{ |c| c | c | c | c |c|c|c|c|c|c|c|c|}
\hline
 & $e_1$& $e_2$&$e_3$ &$e_4$ &$e_5$ & $e_6$& $e_7$\\
\hline
$e_1$ & $0$ &$0$&$0$&$0$&$0$&$0$&$-e_1$ \\ 
\hline
$e_2$ & $0$ &$0$&$0$&$e_1$&$0$&$3e_2$& $-2e_2$\\ 
\hline
$e_3$ & $0$ &$0$&$0$&$0$&$e_1$&$e_3$&$-e_3$ \\ 
\hline
$e_4$ & $0$ &$-e_1$&$0$&$0$&$0$&$-3e_4$&$e_4$ \\ 
\hline
$e_5$ &  $0$&$0$&$-e_1$&$0$& $0$&$-e_5$&$0$\\ 
\hline
$e_6$ & $0$ &$-3e_2$&$-e_3$&$3e_4$&$e_5$&$0$&$0$ \\ 
\hline
$e_7$ &  $e_1$&$2e_2$&$e_3$&$-e_4$&$0$&$0$&$0$ \\ 
\hline
\end{tabular}
\end{center}

\clearpage

\item $\mathfrak{s}_{7, 3}$, nilradical $\mathfrak{n}_{6, 18}$,  where $\mathfrak{s}_{7, 3, A, B} \cong \mathfrak{s}_{7, 3, A', B'}$  iff $AB'=A'B$

\begin{scriptsize}
\begin{center}
\begin{tabular}{ |c| c | c | c | c |c|c|c|c|c|c|c|c|}
\hline
 & $e_1$& $e_2$&$e_3$ &$e_4$ &$e_5$ & $e_6$& $e_7$\\
\hline
$e_1$ &$0$&$-e_3$&$-2e_4$&$3e_5$&$0$&$0$&$(B-2A)e_1$ \\ 
\hline
$e_2$ &$e_3$&$0$&$0$&$0$&$e_6$&$0$&$(3A-2B)e_2$\\ 
\hline
$e_3$ &$2e_4$&$0$&$0$&$3e_6$&$0$&$0$&$(A-B)e_3$ \\ 
\hline
$e_4$ &$-3e_5$&$0$&$-3e_6$&$0$&$0$&$0$& $-Ae_4$\\ 
\hline
$e_5$ &$0$&$-e_6$&$0$&$0$&$0$&$0$&$(B-3A)e_5$\\ 
\hline
$e_6$ &$0$&$0$&$0$&$0$&$0$&$0$&$-Be_6$ \\ 
\hline
$e_7$ &$(2A-B)e_1$&$(2B-3A)e_2$&$(B-A)e_3$&$Ae_4$&$(3A-B)e_5$&$Be_6$&$0$ \\ 
\hline
\end{tabular}
\end{center}
\end{scriptsize}
\vspace{3mm}

\vspace{2mm}

\item $\mathfrak{s}_{8, 1}$, nilradical $\mathfrak{n}_{6, 18}$

\begin{small}
\begin{center}
\begin{tabular}{ |c| c | c | c | c |c|c|c|c|c|c|c|c|}
\hline
 & $e_1$& $e_2$&$e_3$ &$e_4$ &$e_5$ & $e_6$& $e_7$ & $e_8$\\
\hline
$e_1$ &$0$&$-e_3$&$-2e_4$&$3e_5$&$0$&$0$&$-2e_1$& $e_1$\\ 
\hline
$e_2$ &$e_3$&$0$&$0$&$0$&$e_6$&$0$&$3e_2$&$-2e_2$\\ 
\hline
$e_3$ &$2e_4$&$0$&$0$&$3e_6$&$0$&$0$&$e_3$& $-e_3$\\ 
\hline
$e_4$ &$-3e_5$&$0$&$-3e_6$&$0$&$0$&$0$& $-e_4$&$0$\\ 
\hline
$e_5$ &$0$&$-e_6$&$0$&$0$&$0$&$0$&$-3e_5$&$e_5$\\ 
\hline
$e_6$ &$0$&$0$&$0$&$0$&$0$&$0$&$0$& $-e_6$\\ 
\hline
$e_7$ &$2e_1$&$-3e_2$&$-e_3$&$e_4$&$3e_5$&$0$&$0$ &$0$\\ 
\hline
$e_8$ &$-e_1$&$2e_2$&$e_3$&$0$&$-e_5$&$e_6$&$0$&$0$ \\
\hline
\end{tabular}
\end{center}
\end{small}

\end{itemize}



\subsection{Classification of Levi decomposable algebras}\label{levidecomposable}

The Levi decomposable algebras up to  dimension $9$ were classified by Turkowski \cite{turko, turko2}, 
with some small omissions identified by Campoamor-Strusberg \cite{camp}. A detailed description of the classification, and 
 general properties of Levi decomposable algebras is contained in \v{S}nobl and Winternitz \cite{levi}. 

Recall that a Levi decomposable algebra $\mathfrak{g}$ is the semidirect sum of a semisimple Lie algebra $\mathfrak{p}$, called the {\it Levi factor}, and 
the radical $\mathfrak{r}$:
\begin{equation}
\mathfrak{g}= \mathfrak{p} \inplus \mathfrak{r}, ~\text{where}~[\mathfrak{p}, \mathfrak{p}] = \mathfrak{p}, ~[\mathfrak{r}, \mathfrak{r}] \subset \mathfrak{r},~[\mathfrak{p}, \mathfrak{r}] \subseteq \mathfrak{r}.
\end{equation}
An important observation is that $\ad(\mathfrak{p})|_\mathfrak{r}$ is a representation of 
$\mathfrak{p}$ on $\mathfrak{r}$.

The Levi factor of each Levi decomposable subalgebra within this article is $A_1$. It  has basis $\{h, p_+, p_-\}$ with commutation relations
\begin{equation}
[h, p_\pm]=\pm 2 p_\pm, ~ [p_+, p_-]=h.
\end{equation}
For each natural number $j$, $A_1$ has an irreducible representation of dimension $j$ denoted $\rho_j: A_1 \rightarrow \mathfrak{gl}(V)$, 
given,  relative to a basis $\{v_0, v_1,...,$ $v_{j-1}\}$ of $V$, by
\begin{equation}
\begin{array}{lllll}
\rho_j (h) v_k =(j-1-2k)v_k, ~ \rho_j(p_-)v_k=v_{k+1},\\
 \rho_j(p_+)v_k=k(j-1)v_{k-1}, k\in \{0, 1,...,j-1\}.
\end{array}
\end{equation}
 
Below, we describe just those Levi decomposable algebras that  are indecomposable and appear in this article, using the notation of \cite{levi}. For each such Levi decomposable algebra, we include a full multiplication table, identify its radical, and describe the representation $\rho$ of the Levi factor on its radical. 
Not all  Levi decomposable algebras below appear in \cite{levi}, but in cases where they don't occur, we choose notation consistent with that in \cite{levi}.

\vspace{2mm}

 \begin{itemize}
\item  $A_1\inplus_{} 2\mathfrak{n}_{1,1}$, $\rho=\rho_2$
   

\begin{center}
\begin{tabular}{ |c| c | c | c | c |c|}
\hline
 & $h$ & $p_+$ & $p_-$ & $e_1$& $e_2$\\
\hline
$h$ & $0$ & $2p_+$ & $-2p_-$ & $e_1$& $-e_2$\\ 
\hline
$p_+$ & $-2p_+$ & $0$ & $h$ & $0$& $e_1$\\
\hline
$p_-$ & $2p_-$ & $-h$ & $0$ & $e_2$& $0$\\
\hline
$e_1$ & $-e_1$ & $0$ & $-e_2$ & $0$& $0$\\
\hline
$e_2$ & $e_2$ & $-e_1$ & $0$ & $0$& $0$\\
\hline
\end{tabular}
\end{center}
\vspace{2mm}

\vspace{2mm}


 
 
 \item $A_1\inplus_{} 3\mathfrak{n}_{1,1}$, $\rho=\rho_3$

\begin{center}
\begin{tabular}{ |c| c | c | c | c |c|c|}
\hline
 & $h$ & $p_+$ & $p_-$ & $e_1$& $e_2$&$e_3$\\
\hline
$h$ & $0$ & $2p_+$ & $-2p_-$ & $2e_1$& $0$&$-2e_3$\\ 
\hline
$p_+$ & $-2p_+$ & $0$ & $h$ & $0$& $-2e_1$&$e_2$\\
\hline
$p_-$ & $2p_-$ & $-h$ & $0$ & $2e_2$& $-e_3$&$0$\\
\hline
$e_1$ & $-2e_1$ & $0$ & $-2e_2$ & $0$& $0$&$0$\\
\hline
$e_2$ & $0$ & $2e_1$ & $e_3$ & $0$& $0$&$0$\\
\hline
$e_3$ & $2e_3$ & $-e_2$ & $0$ & $0$& $0$&$0$\\
\hline
\end{tabular}
\end{center}
\clearpage

\item  $A_1\inplus \mathfrak{s}_{3,1,  A=1}$,  $\rho=\rho_2\oplus \rho_1$

\begin{center}
\begin{tabular}{ |c| c | c | c | c |c|c|}
\hline
 & $h$ & $p_+$ & $p_-$ & $e_1$& $e_2$&$f_1$\\
\hline
$h$ & $0$ & $2p_+$ & $-2p_-$ & $e_1$& $-e_2$&$0$\\ 
\hline
$p_+$ & $-2p_+$ & $0$ & $h$ & $0$& $e_1$&$0$\\
\hline
$p_-$ & $2p_-$ & $-h$ & $0$ & $e_2$& $0$&$0$\\
\hline
$e_1$ & $-e_1$ & $0$ & $-e_2$ & $0$& $0$&$e_1$\\
\hline
$e_2$ & $e_2$ & $-e_1$ & $0$ & $0$& $0$&$e_2$\\
\hline
$f_1$ & $0$ & $0$ & $0$ & $-e_1$& $-e_2$&$0$\\
\hline
\end{tabular}
\end{center}


 
\item  $A_1\inplus \mathfrak{s}_{4,3, A=B=1}$,  $\rho=\rho_3\oplus \rho_1$

\begin{small}
\begin{center}
\begin{tabular}{ |c| c | c | c | c |c|c|c|}
\hline
 & $h$ & $p_+$ & $p_-$ & $e_1$& $e_2$&$e_3$&$f_1$\\
\hline
$h$ & $0$ & $2p_+$ & $-2p_-$ & $2e_1$& $0$&$-2e_3$&$0$\\ 
\hline
$p_+$ & $-2p_+$ & $0$ & $h$ & $0$& $2e_1$&$2e_2$&$0$\\
\hline
$p_-$ & $2p_-$ & $-h$ & $0$ & $e_2$& $e_3$&$0$&$0$\\
\hline
$e_1$ & $-2e_1$ & $0$ & $-e_2$ & $0$& $0$&$0$&$e_1$\\
\hline
$e_2$ & $0$ & $-2e_1$ & $-e_3$ & $0$& $0$&$0$&$e_2$\\
\hline
$e_3$ & $2e_3$ & $-2e_2$ & $0$ & $0$& $0$&$0$&$e_3$\\
\hline
$f_1$ & $0$ & $0$ & $0$ & $-e_1$& $-e_2$&$-e_3$&$0$\\
\hline
\end{tabular}
\end{center}
\end{small}

\vspace{2mm}

 
\item   $A_1\inplus \mathfrak{s}_{4,3, A=1, B=2/3}$,  $\rho=\rho_2\oplus \rho_1\oplus \rho_1$

\begin{small}
\begin{center}
\begin{tabular}{ |c| c | c | c | c |c|c|c|}
\hline
 & $h$ & $p_+$ & $p_-$ & $e_1$& $e_2$&$e_3$&$f_1$\\
\hline
$h$ & $0$ & $2p_+$ & $-2p_-$ & $e_1$& $-e_2$&$0$&$0$\\ 
\hline
$p_+$ & $-2p_+$ & $0$ & $h$ & $0$& $e_1$&$0$&$0$\\
\hline
$p_-$ & $2p_-$ & $-h$ & $0$ & $e_2$& $0$&$0$&$0$\\
\hline
$e_1$ & $-e_1$ & $0$ & $-e_2$ & $0$& $0$&$0$&$e_1$\\
\hline
$e_2$ & $e_2$ & $-e_1$ & $0$ & $0$& $0$&$0$&$e_2$\\
\hline
$e_3$ & $0$ & $0$ & $0$ & $0$& $0$&$0$&$\frac{2}{3}e_3$\\
\hline
$f_1$ & $0$ & $0$ & $0$ & $-e_1$& $-e_2$&$-\frac{2}{3}e_3$&$0$\\
\hline
\end{tabular}
\end{center}
\end{small}


\clearpage

 
\item  $A_1\inplus \mathfrak{n}_{3,1}$,  $\rho=\rho_2\oplus \rho_1$

\begin{center}
\begin{tabular}{ |c| c | c | c | c |c|c|}
\hline
 & $h$ & $p_+$ & $p_-$ & $e_1$& $e_2$&$e_3$\\
\hline
$h$ & $0$ & $2p_+$ & $-2p_-$ & $0$& $e_2$&$-e_3$\\ 
\hline
$p_+$ & $-2p_+$ & $0$ & $h$ & $0$& $0$&$e_2$\\
\hline
$p_-$ & $2p_-$ & $-h$ & $0$ & $0$& $e_3$&$0$\\
\hline
$e_1$ & $0$ & $0$ & $0$ & $0$& $0$&$0$\\
\hline
$e_2$ & $-e_2$ & $0$ & $-e_3$ & $0$& $0$&$e_1$\\
\hline
$e_3$ & $e_3$ & $-e_2$ & $0$ & $0$& $-e_1$&$0$\\
\hline
\end{tabular}
\end{center}


 \item $A_1\inplus \mathfrak{s}_{4,8, A=1}$, $\rho=\rho_2\oplus \rho_1\oplus \rho_1$

\begin{center}
\begin{tabular}{ |c| c | c | c | c |c|c|c|}
\hline
 & $h$ & $p_+$ & $p_-$ & $e_1$& $e_2$&$e_3$&$f_1$\\
\hline
$h$ & $0$ & $2p_+$ & $-2p_-$ & $0$& $e_2$&$-e_3$&$0$\\ 
\hline
$p_+$ & $-2p_+$ & $0$ & $h$ & $0$& $0$&$e_2$&$0$\\
\hline
$p_-$ & $2p_-$ & $-h$ & $0$ & $0$& $e_3$&$0$&$0$\\
\hline
$e_1$ & $0$ & $0$ & $0$ & $0$& $0$&$0$&$2e_1$\\
\hline
$e_2$ & $-e_2$ & $0$ & $-e_3$ & $0$& $0$&$e_1$&$e_2$\\
\hline
$e_3$ & $e_3$ & $-e_2$ & $0$ & $0$& $-e_1$&$0$&$e_3$\\
\hline
$f_1$ & $0$ & $0$ & $0$ & $-2e_1$& $-e_2$&$-e_3$&$0$\\
\hline
\end{tabular}
\end{center}

\vspace{2mm}

\item $A_1 \inplus  \mathfrak{n}_{5, 3}$, $\rho= \rho_4 \oplus \rho_1$

\begin{small}
\begin{center}
\begin{tabular}{ |c| c | c | c | c |c|c|c|c|c|c|}
\hline
 & $h$ & $p_+$ & $p_-$ & $e_1$& $e_2$&$e_3$ &$e_4$ &$e_5$\\
\hline
$h$ & $0$ & $2p_+$ & $-2p_-$ & $0$&$-3e_2$ &$-e_3$ &$3e_4$&$e_5$\\ 
\hline
$p_+$ & $-2p_+$ & $0$ & $h$ &$0$ & $-e_3$&$-6e_5$&$0$&$e_4$\\
\hline
$p_-$ & $2p_-$ & $-h$  & $0$&$0$ &$0$ &$-3e_2$&$3e_5$&$-\frac{2}{3}e_3$\\
\hline
$e_1$ &$0$ & $0$ & $-e_5$ & $0$& $0$&$0$&$0$&$0$\\
\hline
$e_2$ & $3e_2$  &$e_3$  &$0$  & $0$&$0$ &$0$&$e_1$&$0$\\
\hline
$e_3$ &  $e_3$& $6e_5$&  $3e_2$&$0$ &$0$ &$0$&$0$&$e_1$\\
\hline
$e_4$ & $-3e_4$ &$0$ & $-3e_5$ &$0$ & $-e_1$&$0$&$0$&$0$\\
\hline
$e_5$ & $-e_5$ &$-e_4$ &$\frac{2}{3}e_3$  &$0$ & $0$&$-e_1$&$0$&$0$\\
\hline
\end{tabular}
\end{center}
\end{small}

\clearpage

\item $A_1 \inplus  \mathfrak{n}_{5, 2}$,   $\rho= 2\rho_2 \oplus \rho_1$

\begin{small}
\begin{center}
\begin{tabular}{ |c| c | c | c | c |c|c|c|c|c|c|}
\hline
 & $h$ & $p_+$ & $p_-$ & $e_1$& $e_2$&$e_3$ &$e_4$ &$e_5$\\
\hline
$h$ & $0$ & $2p_+$ & $-2p_-$ &$e_1$ & $-e_2$& $0$&$-e_4$ &$e_5$\\ 
\hline
$p_+$ & $-2p_+$ & $0$ & $h$ &$0$ & $e_1$&$0$&$e_5$&$0$\\
\hline
$p_-$ & $2p_-$ & $-h$  & $0$&$e_2$ & $0$&$0$&$0$&$e_4$\\
\hline
$e_1$ &$-e_1$  & $0$ &$-e_2$  & $0$& $0$&$0$&$0$&$0$\\
\hline
$e_2$ & $e_2$ & $-e_1$ &$0$  & $0$& $0$&$0$&$0$&$0$\\
\hline
$e_3$ & $0$ &$0$ & $0$ & $0$& $0$&$0$&$e_2$&$e_1$\\
\hline
$e_4$ & $e_4$ &$-e_5$ & $0$ &$0$ & $0$&$-e_2$&$0$&$e_3$\\
\hline
$e_5$ & $-e_5$ &$0$ & $-e_4$ & $0$&$0$ &$-e_1$&$-e_3$&$0$\\
\hline
\end{tabular}
\end{center}
\end{small}

\vspace{2mm}

\item $A_1 \inplus \mathfrak{s}_{6, 168, A=B=0}$, $\rho = 2\rho_1 \oplus \rho_4$

\begin{small}
\begin{center}
\begin{tabular}{ |c| c | c | c | c |c|c|c|c|c|c|}
\hline
 & $h$ & $p_+$ & $p_-$ & $e_1$& $e_2$&$e_3$ &$e_4$ &$e_5$&$e_6$\\
\hline
$h$ & $0$ & $2p_+$ & $-2p_-$ & $0$&$-3e_2$ &$-e_3$ &$3e_4$&$e_5$&$0$\\ 
\hline
$p_+$ & $-2p_+$ & $0$ & $h$ &$0$ & $-e_3$&$-6e_5$&$0$&$e_4$&$0$\\
\hline
$p_-$ & $2p_-$ & $-h$  & $0$&$0$ &$0$ &$-3e_2$&$3e_5$&$-\frac{2}{3}e_3$&$0$\\
\hline
$e_1$ &$0$ & $0$ & $-e_5$ & $0$& $0$&$0$&$0$&$0$&$-2e_1$\\
\hline
$e_2$ & $3e_2$  &$e_3$  &$0$  & $0$&$0$ &$0$&$e_1$&$0$&$-e_2$\\
\hline
$e_3$ &  $e_3$& $6e_5$&  $3e_2$&$0$ &$0$ &$0$&$0$&$e_1$&$-e_3$\\
\hline
$e_4$ & $-3e_4$ &$0$ & $-3e_5$ &$0$ & $-e_1$&$0$&$0$&$0$&$-e_4$\\
\hline
$e_5$ & $-e_5$ &$-e_4$ &$\frac{2}{3}e_3$  &$0$ & $0$&$-e_1$&$0$&$0$&$-e_5$\\
\hline
$e_6$&$0$&$0$&$0$&$2e_1$&$e_2$&$e_3$&$e_4$&$e_5$&$0$\\
\hline
\end{tabular}
\end{center}
\end{small}

\clearpage

\item $A_1 \inplus \mathfrak{s}_{6, 155, A=1}$, $\rho = 2\rho_1 \oplus 2\rho_2$

\begin{small}
\begin{center}
\begin{tabular}{ |c| c | c | c | c |c|c|c|c|c|c|c|}
\hline
 & $h$ & $p_+$ & $p_-$ & $e_1$& $e_2$&$e_3$ &$e_4$ &$e_5$ &$e_6$\\
\hline
$h$ & $0$ & $2p_+$ & $-2p_-$ &$e_1$ & $-e_2$& $0$&$-e_4$ &$e_5$&$0$\\ 
\hline
$p_+$ & $-2p_+$ & $0$ & $h$ &$0$ & $e_1$&$0$&$e_5$&$0$&$0$\\
\hline
$p_-$ & $2p_-$ & $-h$  & $0$&$e_2$ & $0$&$0$&$0$&$e_4$&$0$\\
\hline
$e_1$ &$-e_1$  & $0$ &$-e_2$  & $0$& $0$&$0$&$0$&$0$&$-3e_1$\\
\hline
$e_2$ & $e_2$ & $-e_1$ &$0$  & $0$& $0$&$0$&$0$&$0$&$-3e_2$\\
\hline
$e_3$ & $0$ &$0$ & $0$ & $0$& $0$&$0$&$e_2$&$e_1$&$-2e_3$\\
\hline
$e_4$ & $e_4$ &$-e_5$ & $0$ &$0$ & $0$&$-e_2$&$0$&$e_3$&$-e_4$\\
\hline
$e_5$ & $-e_5$ &$0$ & $-e_4$ & $0$&$0$ &$-e_1$&$-e_3$&$0$&$-e_5$\\
\hline
$e_6$ &$0$&$0$&$0$&$3e_1$&$3e_2$&$2e_3$&$e_4$&$e_5$&$0$\\
\hline
\end{tabular}
\end{center}
\end{small}

\end{itemize}

\end{document}